\documentstyle{amsppt}
\input psfig.sty
\magnification=\magstep1
\TagsOnRight
\hsize=5.1in                                                  
\vsize=7.8in

\topmatter

\title Poincar\'e-Reidemeister metric, Euler structures, and torsion
 \endtitle
\rightheadtext{}
\leftheadtext{}
\author  Michael Farber* and Vladimir Turaev \endauthor
\address
School of Mathematical Sciences,
Tel-Aviv University,
Ramat-Aviv 69978, Israel\endaddress
\address
 Institut de Recherche Math\'ematique Avanc\'ee,
Universit\'e Louis Pasteur - C.N.R.S., 7 rue Ren\'e Descartes, 67084 Strasbourg, France
\endaddress
\email farber\@math.tau.ac.il, turaev\@math.u-strasbg.fr
\endemail
\thanks{\footnotemark"*" Partially supported by a grant from the 
Israel Academy of Sciences and Humanities and by
the Herman Minkowski Center for Geometry} 
\endthanks
\abstract{In this paper we define
a {\it Poincar\'e-Reidemeister scalar product} on the determinant line of
the
cohomology of any flat vector bundle over a closed orientable 
odd-dimensional manifold.
It is a combinatorial \lq\lq torsion-type"  invariant which 
refines the {\it PR-metric} introduced in \cite{Fa} and contains an
additional sign
 or phase 
information. 
We compute the PR-scalar product in terms of the torsions of {\it   Euler
structures}, 
introduced in
\cite{T1}, \cite{T2}.  We show that the sign of our PR-scalar product
is determined by the Stiefel-Whitney classes and the semi-characteristic of
the manifold. As an application, we compute the Ray-Singer analytic torsion
via  the torsions of Euler
 structures. Another application: a computation of the twisted
  semi-characteristic in terms of the Stiefel-Whitney classes.}
\endabstract
\endtopmatter

\def\buildrel#1\over#2{\mathrel{\mathop{\null#2}\limits^{#1}}}
\def\buildrul#1\under#2{\mathrel{\mathop{\null#2}\limits_{#1}}}

\define\C{{\bold C}}
\define\R{{\bold R}}  
       
\define\Z{{\bold Z}}

\define\Hom{\operatorname{Hom}}

\redefine\det{\operatorname{det}}

\define\im{\operatorname{im}} 
\redefine\Im{\operatorname{Im}}   


\def\<{\langle}
\def\>{\rangle}
\def\kk{\bold k}

 \def\card {{\text {card}}}
 \def\det  {{\text {det}}}

 \def\mod {{\text {mod}}}
 \def\Hom {{\text {Hom}}}
 \def\Im {{\text {Im}}}
 \def\dim   {{\text {dim}}}
 
 \def\Eul {{\text {Eul}}}
 \def\Int {{\text {Int}}}
\define\Det{{\frak {Det}}}
\def\Re {{\text {Re}}}


 \heading{\bf \S 1. Introduction }
\endheading

 Let $F$ be a flat real vector bundle over a closed odd-dimensional
smooth manifold $X$.
  Ray and Singer \cite{RS}  used the Laplace operators  and
their zeta-function regularized determinants to define a norm on the
determinant line of the cohomology $\det \,H^\ast (X;F)$. Ray and Singer
 showed that their norm is topologically invariant. They 
conjectured that for bundles with orthogonal  structure group, this norm 
  coincides with the Reidemeiter
norm on $\det \,H^\ast (X;F)$, defined using a piecewise linear
triangulation of $X$ and the classical Reidemeister-Franz   torsion. 
This conjecture was proven by
 J. Cheeger     and  W. M\"uller  in their celebrated papers \cite{C} and
\cite{Mu}. 

Although the topologically 
invariant Ray-Singer norm is defined for an {\it arbitrary}
 flat real vector bundle  $F$ over $  X$, the combinatorial counterpart, 
the Reidemeister norm, was known only  
for bundles with unimodular structure group. In 1994  W. M\"uller 
\cite{Mu1}
extended the result of  \cite{C}, \cite{Mu} to all unimodular $F$.  

In  \cite{Fa}, it was shown how to construct combinatorially a norm on 
$\det \,H^\ast (X;F)$ for an arbitrary flat real vector bundle $F$ without
the
unimodularity assumption. The construction of \cite{Fa} uses a combination 
of the   Reidemeister torsion with the Poincar\'e duality;
the resulting norm on $\det \,H^\ast (X;F)$ is called the {\it  
Poincar\'e-Reidemeister norm}. It was proven in \cite{Fa}, that this norm
coincides with the Ray-Singer norm for  any $F$.  The
proof  of this theorem  uses fundamental results of J.-M.
Bismut and W. Zhang \cite{BZ}.
 
A different approach to   the Reidemeister torsions   was 
introduced in \cite{T1} - \cite{T3}. It is observed in these papers 
that    the indeterminacy
  of the Reidemeister torsion is controlled by    additional 
structures on the manifold $X$, the {\it  homology orientations }  and the
{\it  Euler structures.} A homology orientation is 
an orientation 
  of the determinant line of real homologies 
$\det \, H_*(X; \R)$.
An Euler structure on $X$ can
be described in terms of an Euler chain on a PL-triangulation of $X$; it
may
also be   described via vector fields    or via 
$Spin^c$-structures (for 3-manifolds), see loc. cit. 
The constructions of \cite{T1}, \cite{T2} yield  
torsions of Euler structures on $X$ which refine the usual  Reidemeister
torsions.

The initial goal of this research was to find a relation between the
 approaches of  \cite{Fa} and 
\cite{T1}, \cite{T2}. 

In this paper we define the {\it    Poincar\'e-Reidemeister scalar
product } on $\det \,H^\ast (X;F)$, which determines the PR-norm 
defined in \cite{Fa} and contains
an additional sign or phase information. We show that the sign of
 the PR-scalar product is determined by the Stiefel-Whitney classes of $F$
and $X$ and the semi-characteristic of $X$.
The   main result of this paper  computes
 the  
Poincar\'e-Reidemeister scalar
product in terms of  the torsions of Euler structures on   $X$. More
precisely,
in the case of even-dimensional $F$, we give a formula expressing the 
PR-scalar product applied to the torsion of an Euler structure $\xi$ on $X$
in terms of a characteristic homology class $c(\xi)\in H_1(X)$ associated
to  
 $\xi$. For odd-dimensional $F$, we establish a similar formula with the
only difference that the torsion of $\xi$ depends also on a choice of a
homology orientation of $X$. Using these formulas 
 and the main result of  \cite{Fa} we
compute the analytic Ray-Singer torsion   in terms of the  Euler
structures. 

As an application, we  compute the residue  mod  2 of the 
twisted semi-characteristic of $X$ with coefficients  in a  flat vector
bundle
with orthogonal structure group. We give a formula for this residue in
terms of the
Stiefel-Whitney classes. (For related formulas, see   \cite{LMP}). 

In order to prove our results we   develop     general algebraic tools,
allowing to treat the sign anomalies, which appear in the formalism of the
determinant lines. In \cite{T1}, the canonical isomorphism between
the determinant lines of a chain complex and its homology  
was modified by introducing an additional sign factor. In this paper we
introduce   more sign   factors in the natural  maps between the
determinant lines  and we show that   these sign  choices are
  compatible. 

The plan of the paper is as follows.
In Section 2 we give the   algebraic material concerning the
determinant lines of the chain complexes and the canonical isomorphisms.
In Section 3 we recall the construction of the Reidemeister torsion
for unimodular  flat vector bundles. We examine the indeterminacy of this
construction and conclude
that there is no indeterminacy if the flat bundle is even dimensional.
In Section 4 we define the
Poincar\'e-Reidemeister scalar product. Here we state Theorem 4.4 
computing the sign of  this scalar product; the proof appears in Section 6.
In Section 5 we recall the notion of  an Euler structure following  
\cite{T2}. In Section 6 we define the  torsions of Euler structures. 
The novelty here (compared to \cite{T1}, \cite{T2}) 
is that we view the torsion as an element of the determinant line of
homology and also in a systematic treatment of the torsions of flat vector
bundles.   In Section 6  we also state our main theorems 
(Theorems 6.2 and 6.4) computing
the  Poincar\'e-Reidemeister scalar product in terms of the Euler
structures.  In Section 7 we establish two important properties of the
torsion: multiplicativity and duality. They are used in the next Section 8
to
prove   Theorems 6.2 and 6.4. In Section 9 we  introduce   dual 
cohomological versions of both the   PR-scalar product and the
torsion of Euler structures and reformulate our main theorem in these
terms. In Section  10 we 
compute the
analytic Ray-Singer torsion   in terms of the  Euler
structures. Finally, in Section 11 we consider the twisted
semi-characteristics of
$X$. 

 \heading{\bf \S 2. Determinant lines
  of chain complexes}\endheading

In this section  we   recall the    
  canonical isomorphism relating the 
determinant line of a chain complex and the determinant line of its
homology. Our formula (cf. (2-2)) contains a sign refinement, suggested in
\cite{T1}, of the standard formula \cite{M2}. We will introduce also some
sign involving factors in the natural 
commutativity and duality maps between the determinant lines. We will
establish a few technical results concerning
 the compatibility   of  these sign involving choices.

\subheading{2.1. Determinant lines} 
 We shall denote by $\kk$ a fixed ground field of characteristic zero.
The most important special cases are $\kk =\R$ and  $\kk=\C$.

If $V$ is a finite dimensional vector space over $\kk$, {\it the
determinant line of} $V$ is denoted by 
$\det \, V$ and is defined as the top
exterior power of $V$, i.e., $\Lambda^n V$, where $n = \dim \, V$. 
The dual line $\Hom_{\kk} (\Lambda^n V,\kk)$ is denoted by $(\det \, V)^{-1}$.
This notation is justified by the obvious equality 
$\det \, V \otimes (\det \, V)^{-1} =\kk$.
If $V=0$   then by definition $(\det \, V)^{-1}=\det \, V=\kk$. 

For a finite dimensional graded vector space
$V = V_0\oplus V_1\oplus   \dots \oplus V_m $, its {\it determinant
line}  $\det \, V$ is defined as the tensor product 
$$\det \, V = \det \, V_0 \otimes (\det \, V_1)^{-1}\otimes \det \, V_2
\otimes
\dots \otimes (\det \, V_m)^{(-1)^m}.$$

\subheading{2.2. Torsion of a chain complex}  Let $C$ be a finite
dimensional
chain complex  $$0\to C_m\buildrel d\over \to C_{m-1}\buildrel d\over \to
 \dots \buildrel d\over \to C_1\buildrel d\over \to C_0\to 0  $$
over $\kk$.
In the theory of  torsions a crucial role is played by a canonical
isomorphism
$$\varphi_{C}: \det \,  C  \to \det \,  H_\ast(C), \tag2-1$$
where    both $C$
and $H_\ast(C)$ are considered   as graded vector spaces. 
The definition of the
mapping $\varphi_{C}$ is as follows.
 Choose for each $q=0,...,m$  
non-zero elements
$c_q\in \det \, C_q $ and $h_q\in \det \, H_q(C)$. Set 
$c=c_0\otimes c_1^{-1}\otimes c_2\otimes \dots \otimes c_m^{(-1)^m}\in
\det \, C $ and  $h=h_0\otimes h_1^{-1}\otimes h_2\otimes \dots
\otimes h_m^{(-1)^m} \in \det \,H_\ast(C)$,
 where   $-1$ in the exponent 
denotes the dual functional; for example, $c_1^{-1}$ is a 
$\kk$-linear mapping $\det \, C_1 \to \kk$ such that $c_1^{-1}(c_1)=1$.
We define $\varphi_{C}$ by  
$$\varphi_C (c) = (-1)^{N(C)}\, [c:h]\,h, \tag2-2$$
where $N(C)$ is a residue modulo 2 defined  below and $ [c:h]$ is a
nonzero element of $\kk$, defined by $$ [c:h] =  \prod_{q=0}^m
 [d(b_{q+1})\hat h_q b_q/\hat c_q]^{(-1)^{ q+1 }}.\tag2-3$$
Here $b_q$ is a sequence of vectors of $C_q$ whose
 image $d(b_q)$ under the boundary homomorphism $d:C_q\to C_{q-1}$
is a basis of $\Im\, d$; the symbol $\hat h_q$ denotes a sequence of
cycles in $C_q$ such that the wedge product of their homology classes 
equals $h_q$;
the symbol $\hat c_q$ denotes a basis of
  $C_q$ whose wedge product  
equals $c_q$; the number
$ [d(b_{q+1})\hat h_q b_q/\hat  c_q]$ is the determinant of the matrix
 transforming $\hat  c_q$ into the 
basis $d(b_{q+1})\hat h_q b_q$ of $C_q$.
 The residue $N(C)$  is defined by
$$N(C)=\sum_{q=0}^m\alpha_q(C)\beta_q(C) \,(\mod\,2), \tag2-4$$
where
$$\alpha_q(C) = \sum_{j=0}^q \dim \, C_j\,(\mod\,2),
\quad  \beta_q(C) = \sum_{j=0}^q \dim \, H_j(C)\,(\mod\,2).\tag2-5$$
We shall  
deal with chain complexes with zero  Euler characteristic   so that the 
residues
 (2-5) vanish for big 
$q$. 

It is clear that $ [c:h]$ is independent of the choice of
$b_q$'s and also that the isomorphism $\varphi_{C} $ 
is independent of the choice of $h_q$'s and $c_q$'s.

Formula (2-2) involves the  sign
 refinement of the standard formula suggested in \cite{T1}.
In the next subsections we   introduce similar signs
 in other natural maps arising in this setting. 
We shall show that these signs   are   compatible with isomorphism (2-1) 
and with each other.
For more information on torsions, see \cite{M2}, \cite{BGS}, and \cite{Fr}.

\subheading{2.3.  The fusion homomorphism} 
For two  finite-dimensional graded vector spaces $V= V_0\oplus
V_1\oplus   \dots \oplus V_m$
 and $W= W_0\oplus W_1\oplus   \dots \oplus W_m$, we
define a canonical isomorphism  
$$\mu_{V,W}:  \det \, V\otimes \det \, W \to    \det  (V\oplus
W), \tag2-6$$
  by  
$$\mu_{V,W}= (-1)^{M(V,W)}  \bigotimes_q \mu_q^{(-1)^q}, \tag2-7$$
where 
$$\mu_q^{+1}=\mu_q: \det \, V_q\otimes \det \, W_q \to    \det \,
(V_q\oplus W_q)$$  is the isomorphism defined by
  $$(v_1\wedge v_2\wedge \dots \wedge v_k )\otimes
(w_1\wedge w_2\wedge \dots \wedge w_l) \mapsto v_1\wedge v_2\wedge
\dots \wedge v_l \wedge w_1\wedge w_2\wedge \dots \wedge w_k,$$
with $k=\dim \, V_q, l=\dim \, W_q  $,  the isomorphism  
$$\mu_q^{-1}:(\det \, V_q)^{-1}\otimes (\det \, W_q)^{-1} \to    \det \,
(V_q\oplus W_q)^{-1}$$ is the transpose of the inverse  of $\mu_q$,
$${M(V,W)}= \sum_{q=1}^m \alpha_{q-1}(V)\,\alpha_q(W) \,\, \in 
\Z/2\Z, \tag2-8$$ 
with
 $$\alpha_q(V) = \sum_{j=0}^q \dim \, V_j\,\,(\mod\,2) \in  \Z/2\Z,\,\,\,\,
 q=0,1,...,m,$$ and $\alpha_q(W)\in  \Z/2\Z $ defined similarly. 

We will call (2-6) {\it the fusion homomorphism.}
 
\proclaim{2.4. Lemma}  
Let $C=(0\to C_m \to C_{m-1}  \to
 \dots  \to C_0\to 0)$ and
$C'=(0\to C'_m \to C'_{m-1}  \to
 \dots  \to C'_0\to 0)$ be  two finite dimensional chain complexes over
$\kk$. Then the
following diagram, involving the canonical isomorphisms (2-1) and (2-6),
is commutative: $$ \CD \det \, C \otimes \det \, C' @>{\varphi_{C}\otimes
\varphi_{C'}}>> \det \, H_\ast(C)\otimes \det \, H_\ast(C')\\
@V{\mu_{C,C'}}VV  @VV{\mu_{H_\ast(C), H_\ast(C')}}V\\ \det \,(C\oplus
C')@>{\varphi_{C\oplus C'}}>>  \det \, H_\ast(C\oplus C') = \det \,
(H_\ast(C)\oplus H_\ast(C')).  \endCD\tag2-9$$
\endproclaim

\demo{Proof} Fix   non-zero 
$c_q\in \det \, C_q $, $h_q\in \det \, H_q $ and $c'_q\in \det
\,C'_q$, $h'_q\in \det \,H'_q$, where $H_q=H_q(C)$  and
$H'_q=H_q(C')$. We obtain non-zero elements $$c=c_0\otimes
c_1^{-1}\otimes  \dots \otimes c_m^{(-1)^m}\in \det \, C
 \,\,\, {\text {and}} \,\,\,  h=h_0\otimes h_1^{-1}\otimes  
\dots \otimes h_m^{(-1)^m}\in \det \, H_\ast, $$   and similarly $c'\in
\det \,C'$ and $h'\in \det \, H'_\ast$. Set $$cc'=(c_0\wedge
c'_0)\otimes (c_1\wedge c'_1)^{-1}\otimes   \dots
\otimes  (c_m\wedge c'_m)^{(-1)^m} \in \det
\,(C\oplus C') $$ and   $$hh'=(h_0\wedge h'_0)\otimes
(h_1\wedge h'_1)^{-1} 
\otimes \dots \otimes 
(h_m\wedge h'_m)^{(-1)^m} \in \det
\,(H_\ast \oplus H'_\ast).$$  
According to   definitions, 
$$ \mu_{H_\ast , H'_\ast } ((\varphi_{C}\otimes \varphi_{C'})
(c\otimes c'))=(-1)^{N(C) +N(C')}
[c:h] \,[c':h']\, \mu_{H_\ast(C), H\ast(C')}(h\otimes h')$$ $$=(-1)^{N(C)
+N(C')+M(H_\ast , H'_\ast)}\, [c:h] \, [c':h']\, hh'.$$ Similarly,
$$(\varphi_{C\oplus C'}  \mu_{C,C'})(c\otimes c') =
(-1)^{M(C ,C')}\varphi_{C\oplus C'} (cc')=
(-1)^{N(C\oplus C')+M(C ,C')} [cc':hh'] hh'.$$
To prove the lemma we should show that
$$(-1)^{N(C) +N(C')+M(H_\ast , H'_\ast)}
[c:h]  [c':h']=(-1)^{N(C\oplus C')+M(C ,C')} [cc':hh'].\tag2-10$$

Let $b_q$ be a sequence of vectors of $C_q$ whose
 image $d(b_q)$ under the boundary homomorphism $d:C_q\to C_{q-1}$
is a basis of $\Im\, d$. Similarly choose a sequence $b'_q\subset C'_q$
for each $q$.
By definition,
$$\frac{[cc':hh']}{[c:h][c':h']}=
\prod_{q=0}^m \left ( \frac{[(db_{q+1}db'_{q+1}\hat h_q \hat h'_q
b_qb'_q)/c_q]}
{[(db_{q+1}\hat h_q b_q)/c_q]\, [(db'_{q+1}\hat h'_q b'_q)/c_q]}
\right )^{(-1)^{ q+1 }}.$$
The $q$-th factor on the right-hand side is equal to 
 $$(-1)^{ \card (b'_{q+1}) \dim\,H_q +\card (b_{q}) \dim\,H'_q
+\card (b_{q}) \card (b'_{q+1})}.\tag2-11$$ 
  Since  $$\card (b_{q+1}) \equiv \alpha_q(C) +\beta_q(C)\,(\mod\,2),
\quad  \dim\,H_q \equiv
\beta_q(C) + \beta_{q-1}(C)\,(\mod\,2)$$ and similarly for $C'$, we obtain
that the
product of the signs (2-11) equals $(-1)^y$, where
$$y = \beta_{m}(C)\beta_m(C')+$$
$$+\sum_{q=0}^m
\{\beta_q(C)\alpha_q(C')+\beta_{q-1}(C)\beta_q(C')+
\alpha_{q-1}(C)\alpha_q(C') +\alpha_{q-1}(C)\beta_{q-1}(C')\}. $$
It is easy to check that 
$$y\equiv N(C)+N(C') + M(C,C')-N(C\oplus C')-  M(H_\ast(C),H_\ast(C'))
\, (\mod\, 2).$$
This implies (2-10) and the lemma. 
 \qed

\enddemo

\subheading{2.5. Duality operator $D$} 
Let 
$V = V_0\oplus V_1\oplus   \dots \oplus V_m $
be a finite dimensional graded vector space over $\kk$ with odd   $m$.
We define the  {\it  dual graded vector space} over $\kk$
by $V'=  V'_0\oplus V'_1\oplus   \dots \oplus V'_m $
where
$V'_q=(V_{m-q})^*=\Hom_{\kk} (V_{m-q},\kk)$.
We define a {\it duality operator} 
$$D=D_V:\det\,V \to \det\,V'$$ 
as follows. Let $v_q\in \det \, V_q$ be a volume element determined 
by a basis of $V_q$ and let $v'_{m-q}\in \det \, V'_{m-q}$ be the volume
element determined by the dual basis of $V'_{m-q}$, for $q=0,1,...,m$. 
Then   $$D(v_0\otimes v_1^{-1}\otimes v_2\otimes \dots \otimes v_m^{-1}
)=(-1)^{s(V)} v'_0\otimes (v'_1)^{-1}\otimes v'_2\otimes \dots \otimes
(v'_m)^{-1},$$ where
the residue $s(V)\in \Z/2\Z$ is given by
$$s(V)=\sum_{q=1}^m\alpha_{q-1} (V)\, \alpha_q (V)
+\sum_{q=0}^{(m-1)/2}\alpha_{2q}(V) .$$ 
Recall that
$\alpha_q(V) = \sum_{j=0}^q \dim \, V_j\,\,(\mod\,2)$.
It is easy to check that $D=D_V$ does not depend on the choice of 
$v_q\in \det \, V_q\, (q=0,1,...,m)$.

In the next lemma we shall use the notion of a dual chain complex. 
For a chain complex $C=(0\to C_m \to C_{m-1}  \to
 \dots  \to C_0\to 0)$ over $\kk$ the {\it  dual chain complex}
$C'=(0\to C'_m \to C'_{m-1}  \to
 \dots  \to C'_0\to 0)$ is defined by 
$C'_q=(C_{m-q})^* $. The boundary homomorphism 
$C'_{q+1} \to C'_q$ is defined to be 
$(-1)^{m-q} d^*_{m-q-1}$ where $d_{m-q-1}$ is the boundary homomorphism
$C_{m-q} \to C_{m-q-1}$. For odd $m$, the construction  above
yields a duality operator
$D_C:\det\,C \to \det \,C'$. 

\proclaim{2.6. Lemma}  
Let $C=(0\to C_m \to C_{m-1}  \to
 \dots  \to C_0\to 0)$ be a   finite dimensional chain complex with
odd $m$ and with $\chi(C)\equiv 0\, \mod\, 2$, and let  $C'=(0\to C'_m \to
C'_{m-1}  \to
 \dots  \to C'_0\to 0)$ be  the dual chain complex. Then the
following diagram, involving the canonical isomorphisms   (2-1),
is commutative: 
$$
\CD
\det \, C @>{D_C}>> \det \, C'\\
@V{\varphi_C}VV @VV{\varphi_{C'}}V\\
\det \, H_\ast(C) @>{D_{H_\ast(C)}}>> \det \, H_\ast(C').
\endCD\tag2-12 
$$
\endproclaim

Note that the duality between $C$ and $C'$ induces a duality between
the graded vector spaces $H_\ast(C)$ and $ H_\ast(C')$
so that we can consider the duality operator ${D_{H_\ast(C)}}$. 

\demo{Proof} Lemma  2.6 is a sign-refined 
version of the standard duality  for torsions of chain complexes  (see 
\cite{M1}).  For the computation of   signs, 
see Lemma 7 in the Appendix to \cite{T1}. \qed
\enddemo

\proclaim{2.7. Lemma} 
(1) Let $V= V_0\oplus  \dots \oplus V_m$ and $W= W_0\oplus    \dots
\oplus W_m$ be finite-dimensional graded $\kk$-vector spaces   such that 
$\alpha_m(V)=\alpha_m(W) =0\in \Z/2\Z$.  Then 
the following diagram is commutative:
$$
\CD
\det\,V\otimes \det\,W@>{\mu_{V,W}}>>\det\,(V\oplus W)  \\
   @VSVV     @VV{\det(s)}V\\
 \det\,W\otimes \det\,V @>{\mu_{W,V}}>>\det\,(W\oplus V).
\endCD\tag2-13
$$
Here $s$ denotes the natural map $ V\oplus W \to W\oplus V$
 interchanging the summands and 
$S$ interchanges the factors $v\otimes w \mapsto w\otimes v$.

(2) For $V$ and $W$ as in (1) with odd   $m$,
  the following diagram involving  the dual graded vector spaces $V'$ and
$W'$  
 and  the canonical isomorphisms $D$ and $\mu$ is commutative:
$$
\CD
\det\,V\otimes \det\, W @>{D_{V}\otimes D_{W}}>> \det\,V'\otimes
\det\,W'\\ @V{\mu}VV     @VV{\mu}V\\
\det\,(V\oplus W) @>{D_{V\oplus W}}>> \det\,(V'\oplus W'). 
\endCD\tag2-14
$$
 
(3) For any triple of
 finite-dimensional graded vector spaces $U, V,W$, the diagram
$$
\CD
\det\,U \otimes \det\,V\otimes \det\,W @>{\mu_{U,V}\otimes 1}>>
\det\,(U\oplus V)\otimes \det (W)\\
@V{1\otimes \mu_{V,W}}VV   @VV{\mu_{U\oplus V,W}}V\\
\det\,U \otimes \det\,(V\oplus W)@>{\mu_{U,V\oplus W}}>> \det \,(U\oplus
V\oplus W)
\endCD\tag2-15
$$
is commutative.

\endproclaim
\demo{Proof} Statement (1) is equivalent to
$$M(V,W) +M(W,V) +  \sum_q (\dim\, V_q) ( \dim\, W_q) 
\,\equiv \,0 \, (\mod\, 2),$$
which follows easily.

Statement (2) follows from 
$$M(V',W') \equiv  M(W, V)\, \, (\mod \, 2)$$
(using $\alpha_m(V)=\alpha_m(W) =0$) and then
$$M(V, W)+M(W,V)= s(V\oplus W) +s(V)+s(W)\qquad (\mod \, 2)$$

Statement (3) follows from the easy equality
$$M(U,V) +M(U\oplus V,W) = M(V,W)+ M(U,V\oplus W).  $$
\qed
\enddemo

 \heading{\bf \S 3. The Reidemeister torsion}\endheading 
In this section   we   discuss the classical
 construction of the 
 Reidemeister torsion of a flat unimodular bundle. We view this torsion  as
 an element of the determinant line of the  homology of the bundle. 
 We show that the torsion
   has no indeterminacy in the case of an even-dimensional    bundle  and
 has a sign indeterminacy
 in the case of an odd-dimensional   bundle.

\subheading{3.1. Torsion of a unimodular flat vector bundle} 
Let $F$ be a flat $\kk$-vector bundle   over a finite
connected CW-space $X$.  Recall a definition of the homology of $X$
with coefficients in   $F$. Orient all   cells
of $X$.
For a cell $a$ of $X$, denote by $ \Gamma(a, F)$ the vector space
of flat sections of  $F$ over $a$. (Clearly, $\dim\, \Gamma(a, F)=\dim\,
F$.)
The vector space of $q$-chains in $X$ with values in   $F$
 is defined by $$C_q(X; F)  = \bigoplus_{\dim \, a=q}   \,
\Gamma(a, F) . \tag3-1$$ 
The boundary homomorphism $C_q(X; F)\to C_{q-1}(X; F)$ is defined by
restricting
 the flat sections  to the faces
with the signs determined in the usual way 
by the orientations of the cells.
Denote the resulting chain complex by
$C = C_\ast(X;  F)$ and set 
$H_\ast(X;F)= H_\ast(C)$. The graded vector space $H_\ast(X;F)$ is a
homotopy invariant  of the pair $(X,F)$.

Recall    the      Reidemeister-Franz
construction of the torsion of   $(X,F)$.
 We consider here only the case of unimodular $F$, for the
general case, see Remark 3.4 and  Section 6. 
The bundle $F$ is called unimodular,
  if its top exterior power $\Lambda^{\dim\, F} F$ is  a trivial   flat
vector
bundle.
The bundle $F$ 
is unimodular iff it has a flat volume form, i.e.,
a linear volume form on each fiber $F_x, \,x\in X$
invariant
under the parallel transport 
 along any path  in $X$.  Fix such a form
$\omega$.     For
every cell $a$ of $X$ choose  a basis
of $ \Gamma(a,F)$ of  $\omega$-volume 1. 
The concatenation of
these bases over all $q$-dimensional cells gives a basis  in
$C_q(X;F)$ via (3-1). The wedge
product of the elements of this basis yields  
 a non-zero element
$c_q\in   \det \,C_q(X;F)$. 
Set $$\tau(X;F)=\varphi_{C}
(c_0\otimes c_1^{-1}\otimes c_2\otimes \dots \otimes c_m^{(-1)^m})
\in \det \,  H_{\ast}(X;F)  $$
where $m=\dim\, X$ and $\varphi_{C}$ is the 
  isomorphism
$ \det \,  C  \to \det \,  H_\ast(X;F) $
constructed in Section 2.2. In particular,
if $H_\ast(X;F)=0$, then $\det \,  H_\ast(X;F)=\kk$ and 
 $\tau(X;F)\in \kk$ is the Reidemeister-Franz torsion of the pair
$(X,F)$.

The definition of $\tau(X;F)$ involves certain choices.  
Note first that $\tau(X;F)$
does not depend on the choice of   $\omega$-volume 1 bases in $
\{\Gamma(a,F)\}_a$. If  we replace $\omega$ with $k\omega$ for a
non-zero $k\in \kk$, then 
  the torsion $\tau(X;F)$ is
multiplied by $k^{-\chi(X)}$ where $\chi  $ is the    Euler 
characteristic.    
Another   indeterminacy comes from orders and orientations of
the  cells. To apply (3-1), we need to  order the 
$q$-cells of $X$; a permutation in  this order leads to
multiplication of  $\tau(X;F)$ by $(-1)^{\dim \, F}$. Finally, when we
invert the orientation of a cell of $X$ (used in the definition of the
boundary homomorphisms), the torsion $\tau(X;F)$ is also
multiplied by $(-1)^{\dim \, F}$. 
We
sum up this discussion in the following  lemma.

\proclaim{3.2. Lemma} 
Let $F$ be a unimodular flat  vector bundle   over a finite
connected CW-space $X$ with $\chi(X)=0$. The element
$\tau(X;F)\in \det \,  H_\ast(X;F) $ is well defined up to
multiplication
by $ (-1)^{\dim \, F}$.
In particular, if
$F$ is even-dimensional, then 
$\tau(X;F)$ is a well defined element of  $\det \,  H_\ast(X;F) $.
\endproclaim

A fundamental property of the torsion  is its 
 combinatorial invariance which allows to consider the torsions of flat
vector bundles over  PL-manifolds. We have
the following  version of the combinatorial invariance.

\proclaim{3.3. Lemma} Under the conditions of Lemma 3.2, the 
torsion
$\tau(X;F)  $ with indeterminacy given in
Lemma 3.2  is invariant under cell subdivisions of $X$.
\endproclaim
\demo{Proof} The   standard   arguments imply   the combinatorial
invariance of $\tau(X;F)$ modulo $\pm 1$. This yields the
lemma in the case of odd-dimensional $F$. Let us prove the lemma
for even-dimensional   $F$. (We follow the argument given
in \cite{T1}, Section 3.2.1.) Since
a  cellular subdivision is a simple homotopy equivalence,
it is enough to prove  that  
$\tau(X;F)$ 
is invariant
under simple homotopy equivalences. It is well
known that any simple homotopy equivalence may be presented as a
composition of elementary cellular expansions and contractions. 
Therefore it suffices to consider one such transformation. Assume that a
CW-space $X'$ is obtained from $X$  by attaching a closed
$j$-dimensional ball $D$ along a cellular mapping of a closed
$(j-1)$-dimensional ball $D'\subset \partial D$ into $X$. The 
cellular structure in $X'$ is obtained from the one in  $X$ by adding two
open cells  $ a=\Int D$ and $b=\partial D \backslash D'$. The flat vector
bundle $F$ over $X$ extends to a  flat vector bundle $F'$ over $X'$.
Clearly, $H_\ast(X;F)=
  H_\ast(X';F')$ and we should prove that 
$\tau(X;F) =\tau(X';F')$.

 We orient and numerate the cells of
$X'$ (in each dimension) so that the newly attached cells $a,b$ appear at
the very end. 
Denote the chain complex $C_\ast(X';  F')$ and its subcomplex 
 $C_\ast(X;  F)$
by $C'$ and $C$, respectively.
It is clear that $C'_q=C_q$ for $q\neq j, j-1$ and $C'_j=C_j
\oplus \Gamma(a, F')$, $C'_{j-1}= C_{j-1} \oplus
\Gamma(b, F')$. We choose  a flat volume form on $F'$ and volume 1
bases in       $\Gamma(a, F')$,   $  \Gamma(b, F')$, and $C_q $,
 as in   Section 3.1.  Denote these bases   by
 $A$,   $B$, and $\hat c_q$, respectively. Note that $\card\, A=\card\,
B=\dim \, F$.   Choose for each $q $  a
non-zero element  $h_q\in \det \, H_q(C)=\det \, H_q(C')$. 
Choose  a sequence of vectors $b_q$ in $C_q$ whose
 image   under 
the boundary homomorphism $d_{q-1}:C_q\to C_{q-1}$
is a basis of $\Im\, d_{q-1}$. It is easy to see that the image of 
the boundary homomorphism $d'_{q-1}:C'_q\to C'_{q-1}$ equals to 
$\Im\, d_{q-1}$ for $q\neq j$ and that $d_{j-1} (b_j), d'_{j-1} (A)$ is a
basis 
of  $\Im\, d'_{j-1}$.  Note that the residues $N(C)$, $N(C')$
introduced in Section 2.2 are both equal to 0.
Now, it follows   from   definitions that
$$\frac {\tau(X';F')}{\tau(X;F)} =$$
$$=
\left  (\frac {[d_{j}(b_{j+1})\hat h_j b_j A/\hat c_j A]}{
[d_{j}(b_{j+1})\hat h_j b_j  /\hat c_j  ]}  \right )^{(-1)^j}\times
\left  (\frac {[d_{j-1}(b_{j}) d'_{j-1}(A)\hat h_{j-1}
b_{j-1}/\hat c_{j-1} B]}{ 
[d_{j-1}(b_{j})  \hat h_{j-1}
b_{j-1}/\hat c_{j-1} ]}  \right )^{(-1)^{j-1}}.$$
It is obvious that the first factor on the right-hand side equals 1. 
The second factor on the right-hand side equals
$(-1)^{rs} \varepsilon^r $ where $r=\card\, A=\dim \, F$, $s=\card\, 
\hat h_{j-1} +\card\, 
b_{j-1}$, and $\varepsilon$ is the incidence sign of the  oriented  cells
$a,b$. Since $r$ is even, we obtain ${\tau(X';F')}={\tau(X;F)}$.\qed

 \enddemo

\subheading{3.4. Remark} It is easy  to generalize the definition of  
  $\tau(X;F)  $ to the case of a non-unimodular flat
vector bundle $F$ over a finite
connected CW-space $X$ with $\chi(X)=0$.  This gives  an element 
$\tau(X;F)\in \det \,  H_\ast(X;F) $   defined up to
multiplication
by $ (-1)^{\dim \, F}$ and $\det_F(H_1(X))\subset \kk^*$ where
$\det_F:H_1(X)\to \kk^*$ is the determinant of the
monodromy of $F$. We shall consider a more subtle 
torsion in Section 6. 

  \heading{\bf \S 4. The Poincar\'e-Reidemeister scalar product}
\endheading

In this
  section we introduce the Poincar\'e-Reidemeister
 scalar product on the determinant line 
of the  homology of a 
flat vector bundle over a closed orientable odd-dimensional PL-manifold.
It determines the Poincar\'e-Reidemeister metric, introduced in \cite{Fa},
and carries an additional information in the form of a phase (if $\kk =\C$)
or in the form of a sign (if $\kk=\R$).

\subheading{4.1. The dual flat vector bundle} Let $F$ be a 
flat $\kk$-vector bundle over a  finite connected 
CW-space $X$ with  $\chi(X)=0$. 
Recall  the dual flat vector bundle $F^\ast$. The fiber
of $F^\ast$  over a  point $x\in X$
is the dual  vector space $F^\ast_x=\Hom_{\kk} (F_x,\kk)$. For a path
$\gamma: [0,1]\to X$,  the parallel transport
$  F^\ast_x \to F^\ast_y$
 along $\gamma$ 
is the transpose of the parallel transport $  F_y
\to F_x$   along the inverse path
$\gamma^{-1}$. 

It is clear that for any loop $\gamma$ in $X$ we have
$\det_F(\gamma) \cdot \det_{F^\ast}(\gamma) =1$
and therefore  $F \oplus F^\ast$ is  a 
unimodular   flat vector bundle. Since it is also even-dimensional, we
can  apply the construction of Section   3  to obtain a  well defined 
non-zero element
  $  \tau(X;F\oplus F^\ast)\in \det \, H_\ast(X;F\oplus F^\ast)$.

\subheading{4.2. The  duality operator 
} Let $X$ be a  {\it closed connected oriented
piecewise linear manifold of odd dimension m.} Let $F$ be a
flat $\kk$-vector bundle over  $X$.
The standard  homological intersection pairing 
$$H_q(X;F^*)\otimes H_{m-q}(X;F)\to \kk\tag4-1$$ 
allows us to identify  the dual of $H_{m-q}(X;F)$ with $H_q(X;F^*)$.
Applying
  the construction of Section 2.5 to 
the graded vector space $\oplus_{q=0}^m H_q(X;F)$ we obtain
  a canonical isomorphism $$D: \det \, H_\ast(X;F)\to \det \,
H_\ast(X;F^\ast).\tag4-2$$  By definition,
$D = (-1)^{s(F)} 
\otimes_{q=0}^m \psi_q $ 
where the residue $s(F)\in \Z/2\Z$ is given by
$$s(F)=\sum_{q=0}^m\beta_{q-1}\beta_q
+\sum_{q=0}^{(m-1)/2}\beta_{2q} \,\,\, (\mod\, 2),\qquad
\beta_q=\sum_{q=0}^q\dim \,
H_q(X;F)$$
and $\psi_q$ with even $q$ denotes the isomorphism
$$ \det \, H_q(X;F) \to (\det \, H_{m-q}(X;F^*))^{-1}\tag4-3$$
induced by the intersection form,
while
$\psi_q$  with odd $q$ denotes the isomorphism
$$\psi_q:  (\det \, H_q(X;F))^{-1} \to  \det \,
H_{m-q}(X;F^*)$$
inverse to the transpose of (4-3).  

It is easy to check that $D$ does not depend on the choice of
the
orientation of $X$ and therefore
 can be considered for {\it orientable} manifolds. 
(Hint: $\beta_m\equiv \chi (X)=0 \,(\mod\, 2)$.)   As
an
exercise, the reader may check that  
 $s(F)= s(F^\ast) $ (we shall not use
it).

\subheading{4.3. The Poincar\'e-Reidemeister   pairing}
 Let     $F$ be a
flat $\kk$-vector bundle over a closed connected orientable
odd-dimensional PL-manifold
$X$. Denote by  $\mu $
  the canonical fusion isomorphism
$$ 
   \det \, H_\ast(X; F) \otimes
\det \, H_\ast(X;  F^\ast)\to \det \,
 (H_\ast(X; F)\oplus H_\ast(X; F^\ast)) 
= \det \, H_\ast(X; F\oplus F^\ast)
$$ 
  defined in Section 2.3.
  Consider the bilinear pairing
$$\langle \ ,\ \rangle_{PR}: \det \, H_\ast(X;F)\times \det \, H_\ast(X,
F) \to \kk,\tag4-4$$ given by
$$\langle a ,b \rangle_{PR} =
  \mu (a\otimes D(b)) /  \tau(X;F\oplus F^\ast)  \in \kk,  $$ 
where 
$a,b\in 
\det \, H_\ast(X;F)$ and   $D$ is the  
isomorphism (4-2). In other words, $\langle a ,b \rangle_{PR}$ is an
element of
$\kk$ such that $$\mu (a\otimes D(b))= \langle a ,b \rangle_{PR}\, 
 \tau(X;F\oplus F^\ast).$$ The   pairing (4-4) is called the {\it  
Poincar\'e-Reidemeister scalar product}.  

The Poincar\'e-Reidemeister scalar product determines {\it the
Poincar\'e-Reidemeister metric}
(or norm)
on the determinant line $\det \, H_\ast(X;F)$, which was introduced in
\cite{Fa}.
It is given by
$$a\mapsto \sqrt{|\<a,a\>_{PR}|},\quad a\in \det \, H_\ast(X;F)$$
(the positive square root of the absolute value of $\<a,a\>_{PR}$).
The PR-scalar product contains an additional phase or sign information.

In the sequel we shall compute the Poincar\'e-Reidemeister scalar product
in terms of Euler structures  and their torsions. As an application, we
  describe when this 
  scalar product is positive  definite in terms of the
Stiefel-Whitney classes 
$w_{ 1}(F)\in H^{1}(X,\Z/2\Z)$ and $w_{m-1}(X)\in H^{m-1}(X,\Z/2\Z)$.
Namely, we shall prove the following   theorem.

 \proclaim{4.4. Theorem} Let     $F$ be a 
flat $\R$-vector bundle over a closed connected orientable
PL-manifold
$X$ of odd dimension $m$. If $m \equiv 3 \, (\mod\, 4)$ then the 
Poincar\'e-Reidemeister scalar product on $\det \, H_\ast(X;F)$ is positive  definite.
If $m\equiv 1 \,(\mod\, 4)$, then the Poincar\'e-Reidemeister  
scalar product on $\det \, H_\ast(X;F)$ is positive  definite if  and only if
$$ {\<w_1(F)\cup w_{m-1}(X),[X]\>}\, =  s\chi(X)\cdot \dim \, F\qquad (\mod \, 2), \tag4-5$$
where 
 $s\chi(X)$ is the {\it semi-characteristic} of $X$,  defined by 
$$s\chi(X) = \sum_{i=0}^{(m-1)/2} \dim H_{2i}(X;\R).$$
\endproclaim

Theorem 4.4 implies  that the Poincar\'e-Reidemeister
scalar product
is negative definite if and only if
$m\equiv 1 \,(\mod\, 4)$ and
$$<w_1(F)\cup w_{m-1}(X),[X] >\, =  s\chi(X)\cdot \dim \, F +1 \qquad (\mod \, 2).$$

Theorem 4.4 will be proven in Section 6.

\heading{\bf \S 5. Combinatorial Euler structures}\endheading

In this section we recall     combinatorial Euler structures
on CW-spaces and PL-manifolds following \cite {T2}.

\subheading{5.1. Euler structures on CW-spaces} Let $X$ be a finite
connected CW-space with $\chi(X)=0$.
 An {\it Euler chain in $X$} is
  a 
singular 1-chain $\xi$ in $X$  such that
$$d\xi = \sum_a (-1)^{|a|}p_a \tag5-1$$
where $a$ runs over all cells of $X$ and $p_a $
is a point in   $a$; the symbol $|a|$ denotes the
dimension of $a$.
 The vanishing of the Euler 
characteristic  
guarantees the existence of Euler chains. An {\it Euler structure on
$X$} is an equivalence class of Euler chains with respect to an
equivalence relation which we   now describe.

Suppose that 
$ \xi $ and $ \eta $ are two Euler chains in $X$.
 Additionally to (5-1) we have 
$$d\eta = \sum_a (-1)^{|a|}q_a,\quad \text{where}\quad q_a\in
 a . $$
For each cell $a$ choose a path $\gamma_a$ in $a$ 
joining $p_a$ to $q_a$.
Then the chain
$$\xi -\eta +\sum_a (-1)^{|a|}\gamma_a$$
is a 1-cycle; we   denote by
 $d(  \xi, \eta)$ its homology class in $H_1(X)=H_1(X;\Z) $. 
The class $d(  \xi, \eta) $ is clearly independent of the
choice of the paths $\{\gamma_a\}_a$.
 We   say that the Euler chains $ \xi $ and $  \eta $ are
   equivalent  if   $d(  \xi, \eta)=0 $. 
The set of equivalence classes (i.e., the set of Euler structures on
$X$) is denoted by $\Eul(X )$. Sometimes we shall
 denote an Euler structure
and a representing  it Euler chain by the same letter.

\centerline{\psfig{figure=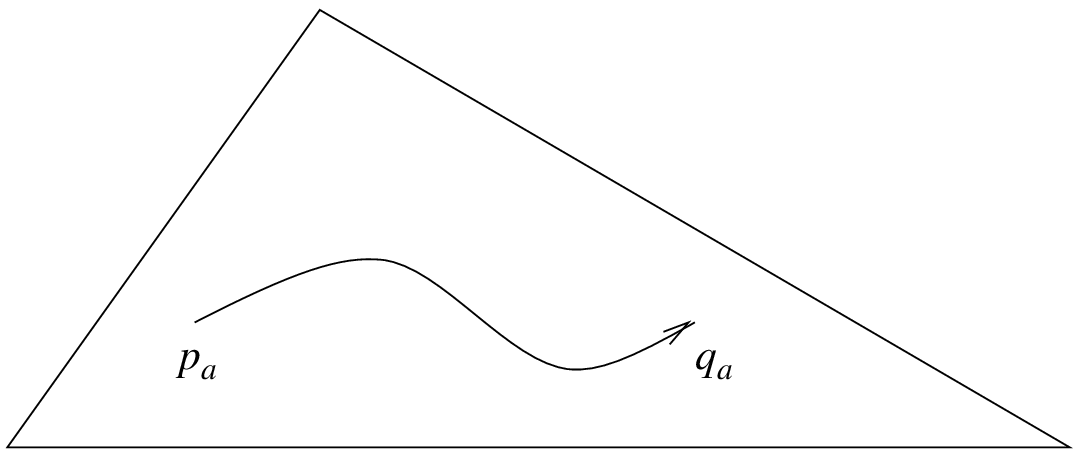,height=1.2in,width=3.0in}}
\vskip 0.5 cm
\centerline{Figure 1. Path $\gamma_a$ }
\vskip 0.5 cm

It is clear that $H_1(X)$ acts on the set $\Eul(X )$: a 1-cycle $h
$ acts on a Euler  chain $\xi$ giving another Euler chain $h+\xi$.
This action of $H_1(X)$ on  $\Eul(X )$ is free and transitive.
We shall use multiplicative notation both for this action and  for the
group operation in $H_1(X)$.

Suppose now that $X'$ is a cellular subdivision of  
$X$.
 Then there is a canonical
bijection
$$\sigma_{X,X'}: \Eul(X)\to \Eul(X'). $$
   It is defined as follows. Let $\xi$ be
 an
Euler chain in $X$ so that (5-1) holds.
Every cell $b$ of $X'$ is contained in a unique cell $a$ of $X$. 
Choose a path $\gamma_b$ in $a$ leading from the point $p_{a}$ to a
certain point in $b$. Set 
$$\xi' = \xi +  \sum_b (-1)^{|b|}\gamma_b $$
where $b$ runs over all cells of $X'$.
It is easy to check that
$\xi'$ is an Euler chain in $X'$.
The correspondence $\xi\mapsto \xi'$
determines a   map 
$\sigma_{X,X'}: \Eul(X)\to \Eul(X')$. 
It is 
$H_1(X)$-equivariant and therefore bijective.

\vskip 1cm
\centerline{\psfig{figure=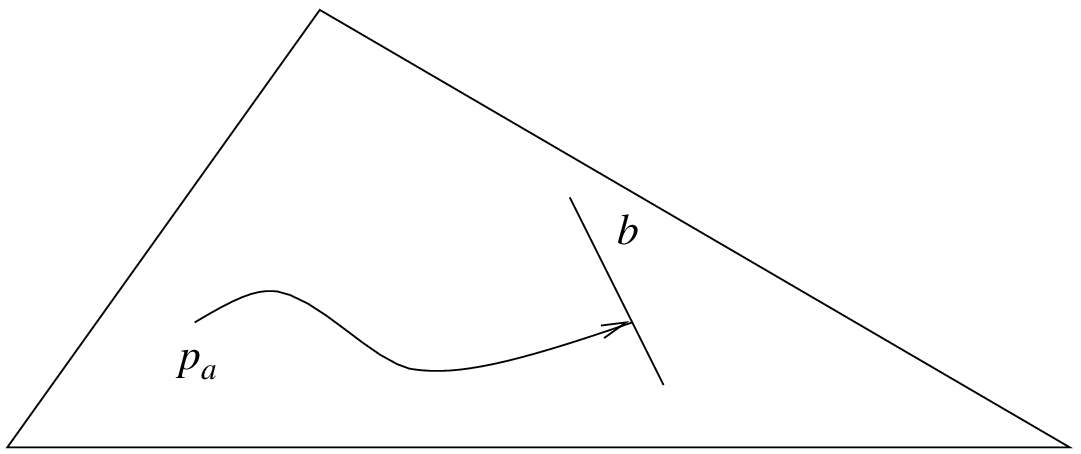,height=1.2in,width=3.0in}}
\vskip 0.5 cm
\centerline{Figure 2. $a, b$ and $\gamma_b$} 
\vskip 0.5 cm
 
\subheading{5.2. Euler structures on PL-manifolds}
Let $X$
be a closed  connected 
PL-manifold  with $\chi(X)=0$.
Each piecewise
linear triangulation $\rho $ of $X$ makes $X$   a
CW-space and allows us to consider the $H_1(X)$-set $\Eul(X,\rho)$.
For a smaller triangulation $\rho'$ we have 
  the equivariant bijection 
$$ \sigma_{\rho,\rho'}: \Eul(X,\rho)\to \Eul(X,\rho') .\tag5-2$$
These sets  and bijections form an inductive system whose 
inductive limit
$$\Eul(X)= \lim_{\rho} \Eul(X,\rho) $$
is   the set of {\it   Euler structures on $X$.}
The group $H_1(X)$ acts on   $\Eul(X)$ freely and transitively.

  For each    Euler structure $\xi$ on 
$X$ we   define its {\it  characteristic class}  $c(\xi)  \in
H_1(X)$ following  \cite{T2}, Section 5.3 and Appendix B. Choose a
PL-triangulation $\rho$ of $X$.  Let $W$ be the   1-chain in $X$
defined by $$W= \sum_{a_0<a_1 \in \rho} (-1)^{|a_0|+|a_1|}
 \langle {\underline {a}}_0,  {\underline {a}}_1 \rangle,$$
where $ a_1$ runs over all simplices of $\rho$, $a_0$ runs over all
proper faces of $a_1$, and 
$\langle {\underline {a}}_0,  {\underline {a}}_1 \rangle$ is a path
in $  a_1$
going from the barycenter ${\underline {a}}_0$
of $a_0$ to the barycenter ${\underline {a}}_1$
of $a_1$. It is easy to check (see \cite{HT})  that 
$$\partial W=(1-(-1)^m) \sum_{a\in \rho } (-1)^{|a | }
  {\underline {a}} $$
where $m=\dim\, X$. Now, any
Euler structure  
 on $ X $ can be presented by   an Euler
chain $\xi$  in
$ (X,\rho)$  such that $ \partial \xi= \sum_{a } (-1)^{|a | }
  {\underline {a}} $. It is clear that $(1-(-1)^m)\, \xi- W$ is a 1-cycle.
Denote its homology class in $H_1(X)$ by $c(\xi)$. 
It follows from \cite {T2}, Lemma B.2.1, that the mapping
$c:\Eul (X,\rho) \to H_1(X)$ commutes with the subdivision isomorphisms
(5-2), i.e., $c\circ \sigma_{\rho,\rho'}=c$.  In this way, we obtain  a
mapping  $c:\Eul (X ) \to H_1(X)$.

Note a few  easy properties of the characteristic class $c $.
 If $m=\dim \,X$ is even, then $c(\xi)$ does not depend on   $\xi$.
If $m$ is odd, then  (in multiplicative notation)
$$c(h\xi) = h^2c(\xi) \tag5-3$$
for any $\xi\in \Eul (X), h\in H_1(X)$. 
For odd $m$, the mod 2 reduction of $c(\xi)$ is independent of 
$\xi$
and equals to the dual of the Stiefel-Whitney class 
$w_{m-1}(X)\in H^{m-1}(X,\Z/2\Z)$. This follows from the
fact that   $W\,(\mod\, 2)$ represents the dual of
$w_{m-1}(X)$, see \cite {HT}.

Using the characteristic class $c$ we   define a mapping
$\xi \mapsto \xi^\ast:\Eul (X)\to \Eul (X)$  by
$$\xi^\ast= (c(\xi))^{-1} \xi.\tag5-4$$
 This mapping is an involution.
It is easy to see this   for odd
$m$. Indeed,  set  $h=c(\xi)$ and observe that
$${\xi^\ast}^\ast=(c(\xi^\ast))^{-1} \xi^\ast=
(c(h^{-1}\xi))^{-1} h^{-1}\xi = 
(h^{-2}h )^{-1} h^{-1}\xi=\xi.$$
For even
$m$, the involutivity of $\ast$ follows from the fact that the 1-cycle
$2W$ is a boundary, see \cite {HT}. 

The involution $\ast$ admits a simple geometric interpretation. Let
$\rho$ be a PL-triangula\-tion of $X$  and let $\rho^\ast$
be the dual cellular decomposition of $X$. Let us represent    $\xi
\in \Eul (X)$ by an Euler chain 
 in
$ (X,\rho)$   denoted by the same letter $\xi$. We can choose this chain
so that  $\partial \xi = \sum_{a\in \rho} (-1)^{|a|}   {\underline {a}}$.
Since the barycenter $  {\underline {a}}$ of $a$ 
 belongs   to
the dual $(m-|a|)$-dimensional cell $a^\ast$,
the 1-chain
$(-1)^m\xi$ is an Euler chain in $ (X,\rho^\ast)$. It   represents
the Euler structure
$\xi^\ast\in \Eul (X)=\Eul (X, \rho^\ast)$  (for a proof, see  \cite
{T2}, Lemma B.2.3).

\centerline{\psfig{figure=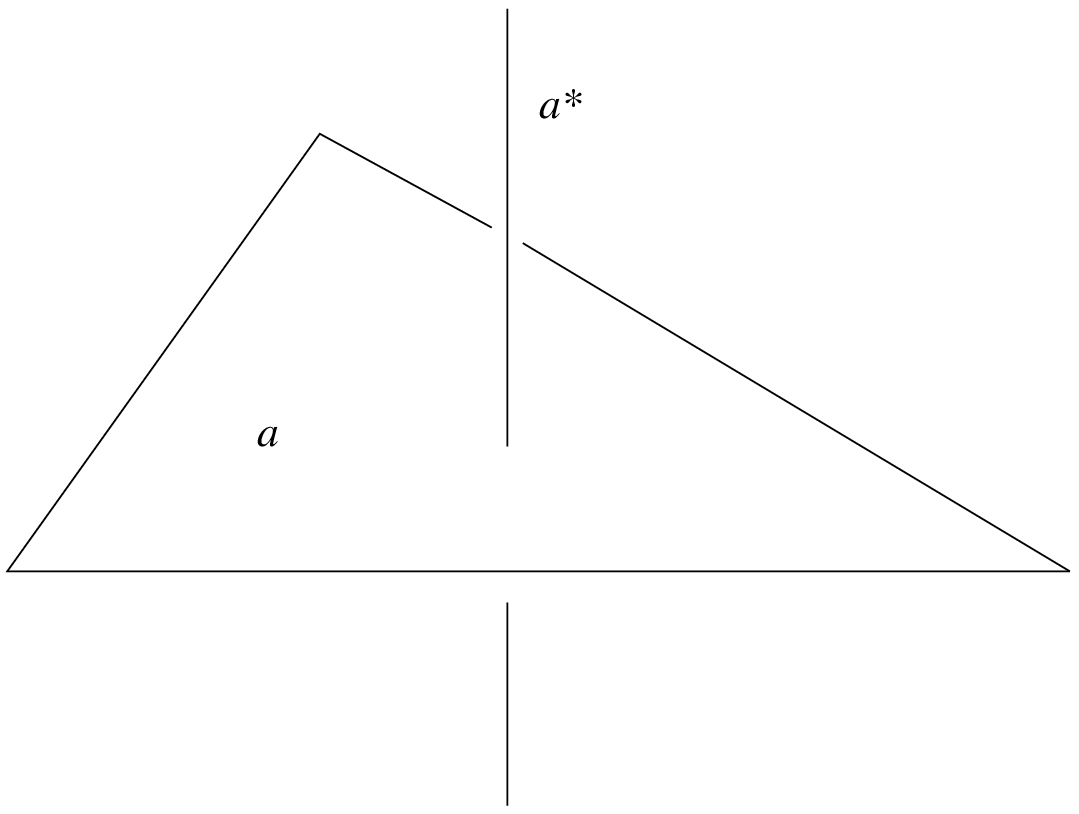,height=1.7in,width=3.0in}}
\centerline{Figure 3.  
Simplex  $a$ and the dual cell $a^\ast$}

 \heading{\bf \S 6. Refined torsions. Main
theorem}\endheading 

In this section we recall     the torsions of Euler structures
  following \cite {T1}, \cite {T2} and state our main theorem
relating them to the 
 Poincar\'e-Reidemeister scalar product. 
We shall    first define the  torsions of Euler structures modulo 
$(-1)^{\dim\, F}$; in particular this gives  well-defined torsions, for 
even-dimensional $F$. For  odd-dimensional  
$F$,
we need to involve additional data (a homology orientation of the base) to
make the torsions of Euler structures 
well-defined.

\subheading{6.1. Torsion of   Euler structures} Let $F$ be a flat vector
bundle   over a finite connected CW-space $X$ with $\chi(X)=0$.
For each Euler structure   $\xi$  on $X$ we define a torsion 
$\tau(X,\xi;F)$ which is an element of the determinant line
$\det \, H_\ast(X;F)$ defined up to multiplication by
$(-1)^{\dim\, F}$. 

As in Section   3  we consider the chain complex 
$C=C_\ast(X;F)$ and the associated torsion isomorphism
$\varphi_{C}: \det \,  C  \to \det \,  H_\ast(X;F) $.
Set $$\tau(X,\xi;F)=\varphi_{C}
(c_0\otimes c_1^{-1}\otimes c_2\otimes \dots \otimes c_m^{(-1)^m})
\in \det \,  H_\ast(X;F)\tag6-1 $$
where $m=\dim\, X$ and   $c_q\in  
\det \,C_q(X;F)\, (q=0,1,...,m)$ are non-zero elements defined as follows.
Fix a point   $x\in X$
and a basis $e_x$ in the fiber $F_x$.
Let
$\beta_a:[0,1]\to X$ be a path connecting $x=\beta_a(0)$ to a point   
$ \beta_a(1) \in   a$. The assumption
$\chi(X)=0$ implies   that the 1-chain
$\sum_{a } (-1)^{|a|}   \beta_a$ (where $a$ runs over all cells of $X$)
is an Euler chain   with boundary $\sum_{a } (-1)^{|a|}  
\beta_a(1)$. We choose the paths $\{\beta_a\}_a$ so that this chain
represents   $\xi$. 
We apply the parallel
transport  to $e_x$ along $\beta_a$ to obtain a
basis in the fiber   $F_{\beta_a(1)}$ and we
extend it   to a basis of flat sections over $a$. The concatenation of
these bases over all $q$-dimensional cells gives a basis  in
$C_q(X;F)$ via (3-1). The wedge
product of the elements of this basis yields  
  $c_q\in   \det \,C_q(X;F)$. 

Let us check the indeterminacy in the definition of $\tau(X, \xi;F)$.   A
different choice of $e_x$
  transforms  the bases in $\{\Gamma(a, F)\}_a$ via one and
the same invertible matrix,   $A$. The torsion
$\tau(X, \xi;F)$   is multiplied by $(\det \, A)^{\chi(X)}=1$
and therefore does not depend on the choice of $e_x$.  
We can replace the path $\beta_a$ by its
composition with a path in $a$ beginning in the point  $\beta_a(1)$.
This does not change the
basis of $\Gamma(a, F)$ constructed   above and therefore does
not change $\tau(X, \xi;F)$. We
can also multiply each $\beta_a$ by   a loop
$\gamma_a: ([0,1], 0,1) \to (X,x,x)$ such that the product 
 $\prod_{a } (-1)^{|a|} \gamma_a $ is homologically trivial.
When  we replace  $\beta_a$ by   $\beta_a\gamma_a$, the element
$c_{|a|} \in   \det \,C_{|a|}(X;F)$  
  is multiplied by $\det_F ([\gamma_a])$ where   $[\gamma_a]\in
H_1(X)$ is the homological class of $\gamma_a$ and 
$\det_F:H_1(X)\to \kk^*$ is the determinant of the monodromy of $F$.
The torsion $\tau(X, \xi;F)$   is multiplied by 
  $ \prod_{a } \det_F(\gamma_a)^{(-1)^{|a|}}=1$ and therefore is not
changed. We can also simultaneously replace the paths
$\{\beta_a\}$ by their compositions $\{\beta_a \gamma\}$ where
$\gamma$ is a path in $X$ leading from a point $y\in X$ to $x$. 
Choosing as $e_y$ the basis in  $F_y$ obtained from $e_x$ by the
parallel transport along $\gamma^{-1}$ we observe that the data 
$y, e_y,   \{\beta_a \gamma\}$ gives rise to the same 
$c_q\in  
\det \,C_q(X;F)\, (q=0,1,...,m)$. Therefore
  $\tau(X,\xi;F)$ does not 
 depend  on  the choice of the point  $x$.
Finally, as in Section 3.1, there is a sign indeterminacy 
$(-1)^{\dim\, F}$ coming from
orders and orientations of the  cells. We conclude that
$\tau(X,\xi;F)$ is  defined up to multiplication by
$(-1)^{\dim\, F}$. In particular,
 for even-dimensional $F$, the torsion $\tau (X,\xi; F)$ is
a well defined element of $\det \, H_\ast(X;F)$.

It follows directly from definitions that 
$$\tau(X,h \xi; F ) = \det_F(h)\cdot \tau(X, \xi;F), $$
 for any $h\in H_1(X)$ and $\xi\in \Eul (X)$.
For unimodular $F$,    we have
$\tau(X, \xi;F)=\tau(X;F)$ where $\tau(X;F)$ is the torsion defined in
Section 3.1.

It follows from \cite {T2},
Lemma 3.2.3
that
the torsion $\tau  (X,\xi;F)$  is invariant under cellular subdivisions
of $X$.  More precisely, if $X'$ is a cellular subdivision of  
$X$ then  
$$\tau(X, \xi;F)=\tau(X', \sigma_{X,X'}(\xi);F)\tag6-2$$
where
 $\sigma_{X,X'}: \Eul(X)\to \Eul(X')$  
is the  canonical
bijection
constructed in Section 5.1. (Note that both parts of (6-2) are
defined up to multiplication by $(-1)^{\dim\, F}$.) This fact allows us to
consider torsions of Euler structures on PL-manifolds.

 \proclaim{6.2.  Main Theorem (even-dimensional case)} 
Let     $F$ be an even-dimensional
flat $\kk$-vector bundle over a closed connected orientable
PL-manifold
$X$ of odd dimension. Then for
  any Euler structure $\xi\in
\Eul(X)$, we have
$$\langle \tau (X,\xi; F) ,\tau (X,\xi; F) \rangle_{PR}
= \det_F(c(\xi)).\tag6-3$$   \endproclaim

Since the mod 2 
reduction of the characteristic class $c(\xi)$ is   dual to the
Stiefel-Whitney  class $w_{m-1}(X)$, Theorem 6.2 implies Theorem 4.4
in the case of even-dimensional $F$.

Using the equality $\xi=c(\xi) \xi^\ast$ we can    reformulate
formula  (6-3)   as follows:
$$\langle \tau (X,\xi; F) ,\tau (X,\xi^\ast; F) \rangle_{PR}=1.$$

To give   similar formulas   for odd-dimensional bundles,  we need   
  a sign-determined version of 
 $\tau (X,\xi; F)$ discussed in the next subsection.

\subheading{6.3. Sign-refined torsion} Let $F$ be an odd-dimensional
flat vector bundle   over a finite connected CW-space $X$ with
$\chi(X)=0$. Assume that $X$ is endowed with an orientation 
$\eta$ of the determinant line of real homologies 
$\det \, H_*(X; \R)$. (Such $X$ is said to be homology oriented.)
Following \cite {T1}, we   introduce for each     $\xi\in \Eul (X)$   a
torsion  $\tau(X,\eta, \xi;F)\in \det \, H_\ast(X;F)$ which has no
indeterminacy.

Let us  orient and  order the cells of $X$. 
Set
$$\tau_0=\varphi_{C}
(c_0\otimes c_1^{-1}\otimes c_2\otimes \dots \otimes c_m^{(-1)^m})
\in \det \,  H_\ast(X;F) $$
where $m=\dim\, X$ and   $c_q\in  
\det \,C_q(X;F)\, (q=0,1,...,m)$ are non-zero elements 
determined by   $\xi$    as in Section
6.1.
Consider the   cellular  chain complex 
 $C_{\R}=C_*(X;\R)$ determined by the trivial line bundle over $X$.
Clearly, $H_*(C)=H_*(X;\R)$.
The     orientation and  order of  the cells of  $X$ 
yield   a basis of  
 $C_{\R}$  which
determines 
  an element   $ c\in  \det \, C_{\R}$.
Recall   the  torsion isomorphism
$ \varphi_{C_{\R}}: \det \, C_{\R}\to \det \, H_*(X;\R)$.
Set $\tau  (X,\eta,\xi ;F)=   \tau_0\in
\det \, H_*(X;F)$ if the   element   $\varphi_{C_{\R}}(c) \in \det \, 
H_*(X;\R)$ defines
  the   orientation $\eta$. In the opposite case set 
$\tau  (X,\eta,\xi ;F)=   -\tau_0\in
\det \, H_*(X;F)$.  It is easy to check that  
  $\tau  (X,\eta,\xi ;F)$ has no
indeterminacy.
 In particular, when we change the orientation or
order of the cells of $X$
 the signs $ (-1)^{\dim\, F}=-1$   appear 
simultaneously in   $ \varphi_{C_{\R}}(c)$ and
$  \tau_0  $    and   cancel each other 
(cf. Section 3.1). 

 Clearly, $\tau(X,  \xi;F) =\pm \tau (X,\eta,\xi ;F)$ 
is the   torsion   discussed in Section 6.1.  
 Note    that  
$\tau  (X,-\eta,\xi ;F)=-\tau  (X,\eta,\xi ;F)$ and 
 $  \tau  (X,\eta,h\xi ;F)=\det_F (h)\,\tau  (X,\eta,\xi ;F)$ for
any $\xi\in \Eul(X)$ and $ h \in H_1(X)$. 

The torsion $\tau  (X,\eta,\xi ;F)$ 
 is invariant under cell subdivisions of $X$, see \cite {T1}, Theorem
3.2.1. (It is to ensure this that we need the signs $ 
(-1)^{N(C) }$ and $  (-1)^{N(C_{\R}) }$  in the definition of the
torsion isomorphisms  $\varphi_{C}, \varphi_{C_{\R}}$.)
 The invariance of
$\tau  (X,\eta,\xi ;F)$  under  cell subdivisions
  allows us to
apply this torsion to PL-manifolds.

\proclaim{6.4.  Main Theorem (odd-dimensional case)} 
Let     $F$ be an odd-dimensional
flat $\kk$-vector bundle over a closed connected orientable
PL-manifold
$X$ of odd dimension $m$. Then for
  any Euler structure $\xi\in
\Eul(X)$ and any homology orientation $\eta$ of $X$, we have
$$\langle \tau (X,\eta,\xi; F) ,\tau (X,\eta, \xi; F) \rangle_{PR}
= (-1)^z \det_F(c(\xi)) \tag6-4$$  
where $z$ is the residue given by   
$$
z=\cases
0,~ {\text  {if}}\,\,\,\dim\,F\,\,\, is \,\,\, even \,\,\, or  \,\,\, 
m\equiv 3\, (\mod\, 4), \\
s\chi(X)\, (\mod\, 2),~ {\text  {if}} \,\,\,
\dim\,F\,\,\, is \,\,\, odd \,\,\, and  \,\,\,  m\equiv 1 \,(\mod\,4). \endcases \tag 6-5
$$ 
\endproclaim

Theorem 6.4 implies the 
identity
$$\langle \tau (X,\eta,\xi; F) ,\tau (X,\eta,\xi^\ast; F)
\rangle_{PR}=(-1)^z.$$

Theorems 6.2 and 6.4 are the main results of this paper. They compute 
the Poincar\'e-Reidemeister scalar product in terms of Euler structures and
their characteristic classes and torsions. 
A proof of Theorems 6.2 and 6.4 is given in Section 8 using the results 
of Section  7.

 \subheading{6.5. Proof of Theorem  4.4}
It follows from Theorems 6.2 and 6.4,
that the Poincar\'e-Reidemeister
scalar product on $\det \, H_\ast(X;F)$ is positive  definite   if and only if
the real number $(-1)^z\det_F(h)$ is positive, where
$h\in H_1(X) $ is a  class whose mod 2 reduction is dual to
$w_{m-1}(X) $ and $z\in \Z/2\Z$ is the residue given by (6-5).
The sign of the non-zero real  number $ \det_F(h)$ is
 equal to $(-1)^{w_1(F) (h)}$
where
 $$w_1(F) (h)={\<w_1(F), w_{m-1}(X)\cap [X]\>}={\<w_1(F)\cup
w_{m-1}(X),[X]\>}. $$
This proves Theorem 4.4 for $m \equiv 1\, (\mod \, 4)$.

It is a theorem of W. Massey \cite{Ma}, Theorem III, that  $w_{m-1}(X)=0$
for any closed orientable smooth manifold $X$ of dimension $m \equiv 3\, (\mod\, 4)$.
This together with the previous argument  gives the claim
of Theorem 4.4 for $m \equiv 3\, (\mod\, 4)$ assuming that $M$ is smoothable. 
Vanishing of the class  $w_{m-1}(X)$ for any closed orientable PL manifold $X$ of 
dimension $m \equiv 3\, (\mod\, 4)$ can be obtained similarly to \cite{Ma}. It also follows
from our arguments used in the proof of Theorem 11.2 (cf. formula (11-2) and Remark 11.4).
This gives our statement for $m \equiv 3\, (\mod\, 4)$ in the PL case.
\qed

\subheading {6.6. Remarks} 1. Any closed oriented manifold
$X$  of odd dimension $m$ has a canonical homology orientation
determined by any basis in $\oplus_{i< m/2} H_i(X; \R)$ followed by
the Poincar\'e dual basis in $\oplus_{i>m/2} H_i(X; \R)$.

2. We could formulate a version of Theorem 6.4   without involving the
sign-refined torsions. Namely, for   odd-dimensional $F$,  we have $$
\langle \tau (X,\xi; F) ,\tau (X,\xi; F) \rangle_{PR}=
(-1)^z \det_F(c(\xi)) \tag6-6 $$
where $z$ is the number defined by (6-5).
This formula makes sense:
although the torsion $\tau (X,\xi; F)$ is
defined up to sign,
the scalar product on the left hand side of (6-6) is   well
defined.
Formula (6-6) directly follows from (6-4).

 \heading{\bf \S 7. Properties of the torsion: multiplicativity and
duality}\endheading

In this section we establish 
 two important properties of the torsion of Euler structures: 
  multiplicativity with respect to direct sums and compatibility
 with  
the duality operator. 
These properties will be used in 
the proof of Theorems 6.2 and 6.4 in Section 8.  

\proclaim{7.1. Theorem}  
Let $F, F'$ be  flat vector bundles  over 
 a finite  connected    CW-space  $X$ with $\chi (X)=0$.   
Let  $$\mu=\mu_{H_\ast(X;F),H_\ast(X;F')}$$
be the canonical fusion isomorphism
$$ 
   \det \, H_\ast(X; F) \otimes
\det \, H_\ast(X;  F')\to \det \, 
(H_\ast(X; F)\oplus H_\ast(X; F')) 
= \det \, H_\ast(X; F\oplus F')
$$ 
  defined in Section 2.3. If both $F$ and $F'$ are even-dimensional then
for any 
    $\xi\in \Eul (X)$
 $$\tau(X,\xi; F\oplus F') =  \mu (\tau(X,\xi; F)\otimes \tau(X,\xi;F')).
\tag7-1$$
If both $F$ and $F'$ are odd-dimensional then
for any  $\xi\in \Eul (X)$ and any homology orientation $\eta$ of $X$,
 $$\tau(X,\xi; F\oplus F') =  \mu (\tau(X,\eta,\xi; F)\otimes
\tau(X,\eta,\xi;F')). \tag7-2$$   \endproclaim

\demo{Proof}  Denote by
$\alpha_q$   the number of cells of $X$ of dimension $\leq q$
and by
${r_q}= \alpha_q-\alpha_{q-1}$   the number of  $q$-dimensional cells
of $X$.   Consider
the  chain complexes  $C=C_*(X; F), C'=C_*(X; F')$, and $\tilde
C=C_*(X; F\oplus F')$. It is clear that $\tilde C=C\oplus C'$.

Let us orient and order the cells of $X$ and fix a spider-like Euler
chain representing $\xi$ as in Section 6.1. The constructions of Section
6.1 provide   bases  in  $C_q$, $C'_q$, and $ \tilde
C_q$.
The   basis  in $C_q $  is formed by a sequence $D_1,...,D_{r_q}$ where
   $D_s$ is a  flat basis of $F$ over the $s$-th
$q$-dimensional cell of $X$.  The basis in  $C'_q $  is formed by a
sequence $D'_1,...,D'_{{r_q}}$ where $D'_s$ is a  flat basis of $F'$ over  
the
$s$-th $q$-dimensional cell of $X$. The   basis in $\tilde C_q$ is formed
by a sequence $D_1,D'_1,D_2, D'_2,...,D_{{r_q}},D'_{{r_q}}$. 
Consider the corresponding wedge products 
$c_q\in \det\,C_q, c'_q\in \det\,C'_q$, and $\tilde c_q\in \det\,\tilde
C_q$. Using  the canonical identification 
$\det\,\tilde
C_q= \det\,C_q \otimes \det\,C'_q$ we obtain 
 $$\tilde c_q=(-1)^{\frac
{(\alpha_q-\alpha_{q-1}-1)(\alpha_q-\alpha_{q-1})} {2}  dd'}
(c_q \otimes c'_q)  $$
where   $d=\dim\, F=\card\,
D_s$  and $  d'=\dim\, F' =\card\, D_s $ for all
$s$.  

Consider   the case where both $d$ and $d'$ are even. 
In this case
$\tilde c_q= c_q \otimes c'_q $ for all $q$. 
By definition,
$$\tau(X,\xi;F)=\varphi_{C}
(c_0\otimes c_1^{-1}\otimes   \dots \otimes c_m^{(-1)^m}), $$
and 
$$\tau(X,\xi;F')=\varphi_{C'}
(c'_0\otimes (c'_1)^{-1}\otimes   \dots \otimes (c'_m)^{(-1)^m}).
$$ 
Lemma  3.3 implies that
$$\mu (\tau(X,\xi;F) \otimes \tau(X,\xi;F'))=$$
$$=
(\varphi_{\tilde C}  \mu_{C,C'} ) 
(c_0\otimes c_1^{-1}\otimes   \dots \otimes c_m^{(-1)^m} 
\otimes  c'_0\otimes (c'_1)^{-1}\otimes   \dots \otimes
(c'_m)^{(-1)^m}) .
$$
By definition of $\mu_{C,C'}$ and by $\tilde c_q= c_q \otimes c'_q $,
 the right-hand side equals 
$$(-1)^{M(C,C')} \varphi_{\tilde C}
(\tilde c_0\otimes (\tilde c_1)^{-1}\otimes   \dots \otimes
(\tilde c_m)^{(-1)^m})=\tau(X,\xi;F\oplus F').$$
Here we use the fact that $\alpha_q(C)=d\cdot \alpha_q$ is even so
that  
  $M(C,C')=0$.

Assume that both $d$ and $d'$ are odd.  
By definition, $\tau(X,\eta, \xi;F)=\varepsilon \tau_0$,
and  $\tau(X,\eta,\xi;F')=\varepsilon \tau'_0$, where 
$$\tau_0=\varepsilon\,\varphi_{C}
(c_0\otimes c_1^{-1}\otimes   \dots \otimes c_m^{(-1)^m}), \,\,\,
\tau'_0=\varepsilon\, \varphi_{C'} (c'_0\otimes (c'_1)^{-1}\otimes  
\dots \otimes (c'_m)^{(-1)^m}),  $$ 
and  $\varepsilon=\pm 1$ is a sign determined by $\eta$ and the
chosen orientations and order of the cells of $X$. It is important that
one and the same sign $\varepsilon$ appears in the expressions for 
$\tau(X,\eta, \xi;F)$ and $\tau(X,\eta,\xi;F')$.
The same argument as above shows that
 $$\mu (\tau(X,\eta,\xi;F) \otimes \tau(X,\eta,\xi;F'))=
\mu (\tau_0 \otimes \tau'_0) 
 =(-1)^{M(C,C')+R} \, \tau(X,\xi;F\oplus F')$$
where
$$R=\sum_{q=0}^m \frac
{(\alpha_q-\alpha_{q-1}-1)(\alpha_q-\alpha_{q-1})} {2}.$$
 It remains to show that 
$M(C,C')+R$ is even. By definition,
$$ M(C,C')=\sum_{q=0}^mdd'\alpha_{q-1} \alpha_{q}\equiv
\sum_{q=0}^m \alpha_{q-1} \alpha_{q} \,(\mod\,2).$$
A direct computation  yields 
$$R\equiv \sum_{q=0}^m \alpha_q-\sum_{q=0}^m \alpha_{q-1} \alpha_{q}
-(\alpha_m+\alpha_m^2)/2\,(\mod\,2).$$
If $P$ (resp. $Q$) is the number of
even-dimensional (resp. odd-dimensional) cells of $X$ then
$P-Q=\chi(X)=0$, $\alpha_m=P+Q=2P$, and 
$\sum_{q=0}^m \alpha_q\equiv  P \equiv  Q \,(\mod\,2)$. This implies
that $M(C,C')+R$ is even and completes the proof of the lemma.\qed

 \enddemo

\proclaim{7.2. Theorem} 
Let $F $ be a flat vector bundle   over   a closed connected orientable PL
manifold $X$ of odd dimension $m$. Let   $\xi\in \Eul (X)$ and let
$ D: \det \, H_\ast(X;F)\to \det \,
H_\ast(X;F^\ast)$ be the   isomorphism  (4-2). If $\dim \, F$ is even then 
$$D(\tau(X,\xi; F )) = \tau (X,\xi^\ast; F^\ast ).\tag7-3$$
If $\dim \, F$ is odd then for any homology orientation $\eta$ of $X$
$$D(\tau(X,\eta,\xi; F )) = (-1)^z \tau (X,\eta,\xi^\ast; F^\ast ),\tag7-4
$$ 
where
$z$ is the number given by (6-5).
\endproclaim

Theorem 7.2 is a  refined 
version of the classical duality  for torsions  due to Franz and Milnor,
see also \cite {T1}, \cite {T2}.

\demo{Proof} Fix an orientation of $X$. Consider first the case of
even-dimensional $F$. Fix a
piecewise linear triangulation $\rho$ of  $X$.  We orient and  order the
simplices of $\rho$ in an arbitrary way.  
Fix a point   $x\in X$.
For each  simplex $a$ of $\rho$, choose
   a path $\beta_a:[0,1]\to X$ connecting $x=\beta_a(0)$ to the
barycenter   of $a$ so   that the 1-chain
$\sum_{a } (-1)^{|a|}   \beta_a$  
represents   $\xi$ in $\Eul (X, \rho)$. 
As in Section 6.1, this chain  and a   basis $e_x$ of the fiber $F_x$ 
determine  an ordered basis of the simplicial chain complex $C =C_*(
(X, \rho);
F)$ and a distinguished element,   $c\in \det\,C$.
By definition, $\tau(X,\xi; F )=\varphi_C(c)$.

To compute the torsion $\tau (X,\xi^\ast; F^\ast )$ we  shall use the 
  dual  
cellular subdivision $\rho^*$ of  $X$. It is
well known that the simplicial chain complex  $C=C_*((X,\rho); F)$ and
the cellular chain complex $C'=C_*((X,\rho^*); F^*)$ are dual to each
other.    Let us provide the cells of $\rho^*$ with the
order and orientation  induced  by the order and orientation of the
simplices of $\rho$.  (To define the induced orientation in the dual cells
we use the orientation of $X$.)
According to the last remark of Section  5.2, the chain
$-\sum_{a } (-1)^{|a|}   \beta_a$ 
represents   $\xi^*$ in $\Eul (X, \rho^\ast)$. 
As in Section 6.1, this chain  and a   basis   of the fiber $F^*_x$ 
determine  an ordered basis of $C'=C_*((X,\rho^*);F^\ast)$
and a distinguished element,   $c'\in \det\,C'$.
By definition, $\tau(X,\xi^*; F^* )=\varphi_{C'}(c')$.
Observe, that if in the role of the basis in 
$F^*_x$  we take the dual basis $e^*_x$, then the basis
in  $C'$ constructed in this way is dual to the basis in $C$ constructed
above. Note that all vector spaces $C'_q$ are even-dimensional
so that
$\alpha_q(C')=0$ for all $q$. Therefore, in
this case $c'=D_C(c)$. It remains to apply  Lemma 2.6 to the complex $C$
(cf. 2.5 and 4.2). This gives 
$$D(\tau(X,\xi; F )) = {D_{H_\ast(C)}} (\varphi_C(c))=\varphi_{C'}
(D_C(c)) =\varphi_{C'}(c')=\tau(X,\xi^*; F^* ).$$

Assume now that $F$ is odd-dimensional. 
As above, we construct
distinguished elements   $c\in \det\,C, c'\in \det\,C'$ and observe that
  $c'=(-1)^{s}  D_C(c)$
  where
  $s\in \Z/2\Z$ is given by
$$s=s(C)=\sum_{q=1}^m\alpha_{q-1}  \alpha_q  
+\sum_{q=0}^{(m-1)/2}\alpha_{2q} \,(\mod\, 2) \tag7-5$$ 
where $\alpha_q$ is the number of simplices of $\rho$ of dimension
$\leq q$. Consider 
the simplicial chain complex $C_{\R}=C_\ast ((X,\rho);\R)$ and the
volume element  $ c_{\R}\in  \det \, C_{\R}$ 
determined by
the     orientation and  order of  the simplices of  $\rho$.
Similarly, consider 
the cellular chain complex $C'_{\R}=C_\ast ((X,\rho^*);\R)$ and the
volume element  $ c'_{\R}\in  \det \, C'_{\R}$ 
determined by
the     orientation and  order of  the cells of  $\rho^*$.
Recall   
the  torsion isomorphisms
 $\varphi_{C_{\R}}:\det \, C_{\R}\to \det \, H_*(X;\R)$ 
 and
 $\varphi_{C'_{\R}}: \det \, C'_{\R}\to \det \, H_*(X;\R)$.
By definition,
$\tau(X,\eta,\xi; F )=\varepsilon \varphi_C(c)$ 
where    $\varepsilon =+1$ if  
the  volume element   $\varphi_{C_{\R}}(c_{\R}) \in \det \, 
H_*(X;\R)$ defines the given homology orientation 
  $\eta$ and $\varepsilon =-1$ otherwise.
 Similarly,
$\tau(X,\eta,\xi^*; F^*)=\varepsilon' \varphi_{C'}(c')$ 
where    $\varepsilon' = +1$ if  
    $\varphi_{C'_{\R}}(c'_{\R}) \in \det \, 
H_*(X;\R)$ defines
  $\eta$  and $\varepsilon' =-1$ otherwise.
As above, $c'_{\R}=(-1)^{s}  D_{C_{\R}}(c_{\R})$
  where $s=s(C_{\R})=s(C)
 \in \Z/2\Z$ is the   residue (7-5).
By Lemma 2.6, 
$$\varphi_{C'_{\R}}(c'_{\R})= \varphi_{C'_{\R}}((-1)^{s } D_{C_{\R}}
(c_{\R}))
= (-1)^{s } D_{H_*(X;\R)}(\varphi_{C_{\R}}(c_{\R})).$$
We can conclude that $ \varepsilon\varepsilon'=(-1)^{s} \nu$
where $\nu=+1$ if the linear mapping
$$D_{H_*(X;\R)}:\det\, {H_*(X;\R)}\to \det\,  {H_*(X;\R)}$$
 preserves the orientation
of the line
$\det\,  {H_*(X;\R)}$ and $\nu=-1$ otherwise.
A computation in 
 \cite{T1}, pp. 178-179 (see also Section 11) shows that $\nu= (-1)^z$
where
$z $ is the number given by (6-5).  
As in the even dimensional case, we
  apply  Lemma 2.6 to the complex $C$
and to the duality operator  ${D_{H_\ast(C)}}=D:\det \,
H_\ast(C) \to \det \, H_\ast(C')$ (cf. 2.5 and 4.2). This gives
$$D(\tau(X,\eta, \xi; F )) = {D_{H_\ast(C)}} (\varepsilon
\varphi_C(c))=  
\varepsilon \varphi_{C'}
(D_C(c))$$ 
$$= \varepsilon \varphi_{C'}((-1)^s c')
=(-1)^{s}   \varepsilon \varepsilon'  \tau(X,\eta,\xi^*;
F^*)= (-1)^z\tau(X,\eta,\xi^*; F^* ).$$
\qed
 \enddemo

 \heading{\bf \S 8. Proof of Theorems 6.2 and 6.4}\endheading

\subheading{8.1. Proof of Theorem  6.2}
Set $T=\tau (X,\xi; F)$. We should prove that
$\langle T ,T \rangle_{PR}  
 =\det_{F}( c(\xi) )$.

Since the bundle
$  F\oplus F^\ast$ is even-dimensional and unimodular,
the torsion  $\tau(X;F\oplus F^\ast) $ is well defined and equals $
\tau(X,\xi; F\oplus F^*)$, for any 
    $\xi\in \Eul (X)$.  By  Theorem 7.1,
 $$\tau(X; F\oplus F^*) = \mu( \tau (X,\xi; F) \otimes \tau(X,\xi;
F^*))= \mu (T\otimes \tau(X,\xi; F^*))  $$ where 
 $\mu $
 is  the canonical fusion isomorphism
$$ 
   \det \, H_\ast(X; F) \otimes
\det \, H_\ast(X;  F^\ast)\to \det \,
 (H_\ast(X; F)\oplus H_\ast(X; F^\ast)) 
= \det \, H_\ast(X; F\oplus F^\ast)
$$ 
  defined in Section 2.3. 
By Theorem 7.2,
$$D(T) = \tau (X,\xi^\ast; F^\ast )=\tau (X, (c(\xi))^{-1}\xi; F^\ast
)$$
$$=\det_{F^*}((c(\xi))^{-1})\,\tau (X,  \xi; F^\ast
)=\det_{F}( c(\xi) )\,\tau (X,  \xi; F^\ast).$$
By definition,
$$\langle T ,T \rangle_{PR} =
  \mu (T\otimes D(T)) /  \tau(X;F\oplus F^\ast)    $$ 
$$= \mu (T\otimes \det_{F}( c(\xi) )\,\tau (X,  \xi; F^\ast)) /
\mu (T\otimes \tau(X,\xi;
F^*)) =\det_{F}( c(\xi) ).$$
\qed

\subheading{8.2. Proof of Theorem  6.4} 
Set $T=\tau (X,\eta, \xi; F)$. 
By  Theorem 7.1,
 $$\tau(X; F\oplus F^*) =  \tau(X,\xi ; F\oplus F^*)=\mu (T\otimes
\tau(X,\eta, \xi; F^*)).  $$
By Theorem 7.2,
$$D(T) = (-1)^z \tau (X,\eta,\xi^\ast; F^\ast )=(-1)^z \det_{F}(
c(\xi) )\,\tau (X,  \xi; F^\ast).$$
Thus,
$$\langle T ,T \rangle_{PR} =
  \mu (T\otimes D(T)) /  \tau(X;F\oplus F^\ast)    $$ 
$$=  \mu (T\otimes (-1)^z \det_{F}( c(\xi) )\,\tau (X, \eta, \xi;
F^\ast)) / \mu (T\otimes \tau(X,\eta,\xi;
F^*)) =(-1)^z \det_{F}( c(\xi) ).$$
\qed

\heading{\bf \S 9. Cohomological torsions and  
the PR-pairing}\endheading

In this section we give   cohomological versions of both the
Poincar\'e-Reidemeister
scalar products and the torsions of Euler structures.
This cohomological formulation is better suited for a comparison with the
analytical approach, see Section 10. 

\subheading{9.1. Cohomology of a flat vector bundle} 
Let $F$ be a flat $\kk$-vector bundle   over a finite
connected CW-space $X$.  Recall a definition of the cohomology of $X$
with coefficients in   $F$. Orient all   cells
of $X$.
As in Section 3.1, for a cell $a$ of $X$, denote by $ \Gamma(a, F)$ the
vector space
of flat sections of  $F$ over $a$.  
The vector space of $q$-cochains in $X$ with values in   $F$
coincides with the vector space of $q$-chains and 
 is defined by $$C^q(X; F)  = \bigoplus_{\dim \, a=q}   \,
\Gamma(a, F) . \tag9-1$$ 
The boundary homomorphism $\delta_q: C^q(X; F)\to C^{q+1}(X; F)$ is defined
as
follows.
  Let $s_a$ be a flat section of $F$ over a 
$q$-cell $a$.
We set 
$$\delta_q (s_a) = \sum_b \varepsilon(a, b) s_a^b $$
where the sum runs over all  $(q+1)$-cells $b$ 
incident to $a$, the sign $\varepsilon(a, b)=\pm 1$ is   determined in the
usual way by the
orientations of $a$ and $b$, and $s_a^b$ denotes the  unique  flat section
over $b$ 
  extending   $s_a$.  (It is understood that each $b$ enters this sum with
multiplicity equal to the number of appearances of $a$ in   $\partial b$.)
 Denote
the resulting cochain complex by $ C^\ast(X;  F)$ and set 
$H^\ast(X;F)= H^\ast(C^\ast(X;  F))$. The graded vector space $H^\ast(X;F)$
is a
homotopy invariant  of the pair $(X,F)$.
 
It is clear   that the 
vector space   $C^q(X;F)$ is dual
to $C_q(X; F^\ast)$, i.e., $C^q(X;F)= \Hom_{\kk} (C_q(X; F^\ast),\kk)$
and the boundary homomorphism  $\delta_q$ introduced above is dual to
the boundary homomorphism $C_{q+1}(X; F^\ast)\to C_q(X; F^\ast)$.
Therefore for each $q$, we have a non-singular evaluation pairing
$$ H^q(X;F)\otimes   H_q(X;F^\ast)\to \kk.$$
These pairings for $q=0,..., \dim\,X$ induce 
a   non-singular  pairing
$$ [\ ,\ ]: \det \, H^\ast(X;F)\otimes \det \, H_\ast(X;F^\ast)\to
\kk.\tag9-2$$
 
\subheading{9.2. Cohomological torsion}  
Let $F$ be a flat $\kk$-vector bundle   over a finite
connected CW-space $X$ with $\chi(X)=0$.  If $\dim\, F$ is odd, then we
additionally
assume that $X$ is provided with a homology orientation (which we suppress
in the
notation). For every $\xi\in \Eul(X)$, we define  the cohomological torsion
$ \tau^\bullet(X,\xi;F)$ as the unique element of $ \det \, H^\ast(X;F) $
such that
$$[\tau^\bullet(X,\xi, F), \tau(X,\xi, F^\ast)] =1 \tag9-3$$
where  
$\tau(X,\xi;F^\ast)\in \det \, H_\ast(X;F^\ast)$  is 
the torsion  defined   in Section 6.
The cohomological torsion
$\tau^\bullet(X,\xi;F)\in \det \, H^\ast(X;F)$ satisfies properties similar
to those of the
homological torsion.
In particular, it is invariant under cell subdivisions and
has no   indeterminacy. 
For any $h\in H_1(X)$, we have
$$\tau^\bullet(X,h\xi;F)= \det_F(h)\cdot \tau^\bullet(X,\xi;F).\tag9-4$$

 \subheading{9.3. Cohomological Poincar\'e-Reidemeister
scalar product} We define a cohomological version of the
Poincar\'e-Reidemeister
scalar product. The norm determined by this scalar product was originally
defined in \cite{Fa}.

Let $F$ be a
flat $\kk$-vector bundle over  a closed connected orientable
  PL manifold $X$
of odd dimension $m$.   
Given $\alpha, \beta\in \det\, H^\ast(X;F)$, we define a number
$ \<\alpha, \beta\>_{PR} \in \kk $
by  
$$\<\alpha, \beta\>_{PR} = \frac{[\alpha,a] \cdot[\beta,b]}{ \<a,
b\>_{PR}}\tag9-5$$ for any nonzero $a,b\in \det\, H_\ast(X;F^\ast)$. Here
 $\<a, b\>_{PR}$ denotes the
homological Poincar\'e - Reidemeister scalar product (defined in Section 4)
and the square brackets denote the     pairing  (9-2).
Formula  (9-5) yields a well-defined bilinear
form on $\det\, H^\ast(X;F)$ called the cohomological
Poincar\'e-Reidemeister
scalar product.

Let us show that the norm, 
$$\alpha\mapsto |\<\alpha,\alpha\>_{PR}|^{1/2},\quad \alpha\in \det\,
H^\ast(X;F),$$
determined by scalar product (9-5), coincides with the Poincar\'e -
Reidemeister norm 
on $\det\, H^\ast(X;F)$,
introduced in \cite{Fa}, section 4.7. This fact will be used in Section 10.

Denote by 
$$\mu^\bullet: \det\, H^\ast(X;F)\otimes \det\, H^\ast(X;F^\ast)\to \det\,
H^\ast(X;F\oplus F^\ast)$$
the canonical isomorphism, defined similarly to (2-6) (ignoring the signs).
Let
$$D^\bullet : \det\, H^\ast(X;F)\to \det\, H^\ast(X;F^\ast),$$
be the Poincar\'e duality isomorphism (a cohomological version of (4-2)). 
In our present notation, the
Poincar\'e-Reidemeister metric on $\det\, H^\ast(X;F)$, which defined
in   Section
4.7 of \cite{Fa},  is given by 
$$\alpha \mapsto
|\mu^\bullet(\alpha\otimes D^\bullet \alpha)/\tau^\bullet(X; F\oplus
F^\ast)|^{1/2},$$ where $\alpha \in \det\, H^\ast(X;F)$. In
order
to prove compatibility with   (9-5) it is enough to show  that for any
$a, b \in \det\, H_\ast(X;F^\ast)$, and $\alpha,\beta\in \det\,
H^\ast(X;F)$ holds
 $$(\mu^\bullet(\alpha\otimes D^\bullet
\beta)/\tau^\bullet(X; F\oplus F^\ast))\cdot  (\mu (a\otimes Db)/\tau(X;
F^\ast\oplus F))= \pm [\alpha, a]\cdot [\beta, b]. $$ 
The left hand side may be rewritten
as
$$\frac {[\mu^\bullet(\alpha\otimes
D^\bullet \beta), \mu (\alpha\otimes D \beta)]}
{[\tau^\bullet(X; F\oplus F^\ast), \tau(X;
F^\ast\oplus F)]}
$$
where the   brackets $[\, , \,]$ denote the pairing (9-2)
for the flat vector bundle 
$ F\oplus F^\ast$.
By (9-3), the denominator of the last expression is equal to $\pm 1$.
It remains to check that  $$[\mu^\bullet(\alpha\otimes D^\bullet \beta),
\mu
(\alpha\otimes D \beta)] =  \pm [\alpha, a]\cdot [\beta, b]$$ where the
brackets $[\ ,\
]$ denote the pairings   (9-2) for 
$F, F^\ast$, and $ F\oplus F^\ast$.
The last equality  follows from
$$[\mu^\bullet(\alpha\otimes D^\bullet \beta), \mu (\alpha\otimes D \beta)]
=  \pm
[\alpha, a]\cdot [D^\bullet\beta, Db] = \pm [\alpha, a]\cdot [\beta, b]$$
(cf. \cite{Fa},
Section 3.4).

\proclaim{9.4.  Main Theorem (cohomological version)} 
Let     $F$ be a 
flat $\kk$-vector bundle over a closed connected orientable 
PL-manifold
$X$ of odd dimension $m$. If $\dim\, F$ is odd, then we additionally assume
that $X$ is provided with a homology orientation.
 Then for
  any   $\xi\in
\Eul(X)$,
$$\langle \tau^\bullet (X,\xi; F) ,\tau^\bullet (X, \xi; F) \rangle_{PR}
= (-1)^z \det_F(c(\xi)), \tag9-6$$  
where $\langle \, ,\, \rangle_{PR}$ is the 
cohomological Poincar\'e-Reidemeister
scalar product and 
$z$ is the
number given by   (6-5). \endproclaim

\demo{Proof} Applying  (9-5)  to  
$$\alpha=\beta=\tau^\bullet(X,\xi;F) \in \det\, H^\ast(X;F)$$ and
$$a=b=\tau(X,\xi;F^\ast)
\in \det\, H_\ast(X;F^\ast)$$ we obtain   $$
 \<\tau^\bullet(X,\xi;F), \tau^\bullet(X,\xi;F)\>_{PR} =$$
$$   
 =\<\tau(X,\xi;F^\ast), \tau(X,\xi;F^\ast)\>_{PR}^{-1}\,\cdot\,
[\tau^\bullet(X,\xi, F),
\tau(X,\xi, F^\ast)]^2. $$
By (9-3),  $[\tau^\bullet(X,\xi, F),
\tau(X,\xi, F^\ast)]=1$. By Theorems 6.2 and 6.4, 
$$\<\tau(X,\xi;F^\ast), \tau(X,\xi;F^\ast)\>_{PR}^{-1}  =(-1)^z
\det_{F^\ast}(c(\xi))^{-1} = 
 (-1)^z \det_{F}(c(\xi)).
$$
This implies the claim of the theorem. \qed
 \enddemo

  \heading{\bf \S 10. Analytic torsion via Euler structures}\endheading

In this section we describe a   
relationship  between the analytic torsion of Ray and Singer
 \cite{RS} and
the combinatorial torsion of   Euler structures. 
The analytic torsion  of a flat vector bundle 
$F$ over a closed odd-dimensional manifold $X$ can be viewed 
  as a norm ({\it the Ray-Singer norm})
 on the determinant line $\det \, H^\ast(X; F)$.
 The main result of this section expresses 
the Ray-Singer norm of the cohomological torsion  
$\tau^\bullet(X, \xi; F) $ of any Euler
 structure $\xi\in Eul(X)$
in terms of the 
monodromy of $F$ along the characteristic class $c(\xi)$. 

\subheading{10.1. Ray-Singer norm}
 We recall the construction of the Ray-Singer norm.
Let $X$ be a closed smooth manifold, 
and let $F$ be a flat real vector bundle over $X$. (Here 
  the ground field $\kk$ is $\R$.)
Choose an arbitrary Riemannian metric 
on $X$ and a smooth metric on $F$. Then the space 
$\Omega^\ast(X;F)$ of
differential forms on $X$ with values in $F$ has a scalar product. 
The flat structure on $F$ determines a flat connection
 $\nabla: \Omega^q(X;F)\to \Omega^{q+1}(X;F)$, so that $\nabla^2=0$. 
We have
$$H^\ast(X;F)=\ker (\nabla)/\im (\nabla)$$
 (the cohomology of the twisted de Rham complex).
Using the Hodge decomposition, 
the cohomology can be embedded into $\Omega^\ast(X;F)$ as
 the space of harmonic forms; this embedding induces a 
norm   $|\cdot|^{RS}$ on the determinant line $\det \, H^\ast(X;F)$. 
The Ray-Singer norm $||\cdot ||^{RS}$ on   $\det \, H^\ast(X;F)$ is
defined by
$$||\cdot ||^{RS} = |\cdot |^{RS} \,
\prod_{q=0}^{\dim X} (\Det \,\Delta'_q)^{(-1)^{q}\, q/2},\tag10-1$$
where $\Det \, \Delta'_q$ denotes the zeta-function regularized determinant
of the Laplacian $\Delta'_q$ acting on the space of $q$-forms   orthogonal
to the harmonic
forms. 
Recall the definition of    $\Det \,\Delta'_q$ following \cite{RS}.
Consider the positive
eigenvalues of the Laplacian $\Delta_q:\Omega^q(X;F)\to \Omega^q(X;F)$
$$0<\lambda_1\le \lambda_2
\le \dots \le \lambda_k\le\dots,\qquad \lambda_k\to \infty$$
and form  the $\zeta$-function
$$\zeta_q(s) =
 \sum_{k=1}^\infty \lambda_k^{-s},\qquad\text{$\Re(s)$ is large}.$$
It is a meromorphic function  holomorphic at $s=0$. Now,
$$ \Det \, \Delta'_q = \exp(-\frac{d}{ds} \zeta_q(s)|_{s=0}).$$

The fundamental property of the Ray-Singer norm (10-1) 
for {\it odd-dimensional}   $X$   is its 
  topological  invariance: it
 does not depend on the choice of   metrics on $X$ and   $F$,  
 used in the
construction. For even-dimensional $X$  this is not the case,
see \cite{BZ} for a detailed description of
the dependence  
of the Ray-Singer norm   on the metrics.

\proclaim{10.2. Theorem (Analytic torsion and Euler structures)} Let 
$X$ be a closed connected orientable smooth manifold of odd dimension  
and let $F$ be a flat $\R$-vector bundle over $X$.
If $\dim\, F$ is odd, then we additionally assume
that $X$ is provided with a homology orientation.
For any Euler structure $\xi\in Eul(X)$,
 the Ray-Singer norm of   its cohomological torsion (cf. 9.2)
$\tau^\bullet(X, \xi; F)\in \det \, H^\ast(X; F)$
 is equal to
the  positive  square root of the absolute value of the monodromy of $F$
along the characteristic 
class $c(\xi)\in H_1(X)$:
\endproclaim
$$||\tau^\bullet(X, \xi;F)||^{RS} = |\det_Fc(\xi)|^{1/2}.\tag10-2$$

In the special case, where the flat bundle $F$ is acyclic, i.e., 
 $H^\ast(X;F)=0$, the torsion 
$\tau^\bullet(X, \xi; F)$ is a real number 
 and Theorem 10.2
yields
$$\prod_{q=0}^{\dim X} (\Det \, \Delta'_q)^{(-1)^{q+1}\, q}= 
\frac{(\tau^\bullet(X, \xi; F))^2}{|\det_Fc(\xi)|}.\tag10-3$$
Note that the RHS of this formula  does not depend on the choice of
 $\xi$, 
this follows directly from (6-3) and   the properties of the torsion.
 
Theorem 10.2 generalizes the classical Cheeger-M\"uller theorem \cite{C},
\cite{Mu1}  concerning the orthogonal flat real bundles $F$  and 
 the (more general) theorem of M\"uller \cite{Mu2}
 concerning the unimodular flat real bundles $F$. 
Note that if $F$ is unimodular then   
$|\det_Fc(\xi)|=1$ and the torsion $\tau^\bullet(X, \xi;F)$ does not
 depend on the choice
of $\xi$.

\demo{Proof of Theorem 10.2}  The main theorem of
 \cite{Fa}, Theorem 3.2, states that  
 the norm on  
 $\det \, H^\ast(X;F)$ associated with the Poincar\'e-Reidemeister 
scalar product 
coincides with the Ray-Singer norm. More precisely,  
 for any $\alpha\in \det \, H^\ast(X;F)$,
$$\<\alpha,\alpha\>_{PR} = \pm \,(||\alpha||^{RS})^2.\tag10-4$$
Substituting here $\alpha = \tau^\bullet(X,\xi,F)$ and using (9-6) we
 obtain  (10-2).\qed
\enddemo

Note that the sign  in (10-4)  is completely described in Theorem 4.4.

\subheading{10.3. Questions} Formula (10-2) 
computes   $\pm \tau^\bullet(X,\xi; F)$ 
 in analytical terms. 
Is there a way to compute  
$\tau^\bullet(X,\xi; F)$ (without the sign indeterminacy)
  using the analytic tools?
 One may expect that the $\eta$-invariant of Atiyah, Patodi and Singer 
will be relevant for
this purpose. 

There is a similar question. 
Suppose that $F$ is a  flat complex bundle.
Then we have the {\it complex torsion}
$\tau^\bullet(X,\xi; F)$ lying in the complex determinant 
line $\det \, H^\ast(X; F)$.
Now, one may also consider $F$ as the real flat bundle
 $F_\R$ and consider the {\it real
torsion} $\tau^\bullet(X,\xi; F_\R)\in \det \, H^\ast(X; F_\R)$.
It can be shown that the real torsion $\tau^\bullet(X,\xi; F_\R)$ 
may be considered as an {\it \lq\lq absolute value"} of the 
complex torsion $\tau^\bullet(X,\xi; F)$; it can be expressed
 in terms of the analytic
torsion of Ray and Singer and the information contained 
in the characteristic class $c(\xi)$, 
using our Theorem 10.2. One may ask how to recover
 {\it the \lq\lq  phase information"}
of the complex torsion $\tau^\bullet(X,\xi; F)$ using the 
analytic tools?

\heading{\bf \S 11. Semi-characteristics of manifolds}\endheading

In this section we will apply the   results obtained above to compute
  the residue mod $\,2$ of the twisted semi-characterictic
of a closed orientable smooth manifold of dimension $\equiv 1\, (\mod\,4)$.

\subheading{11.1. Twisted semi-characteristics}
 Let $F $ be a flat     vector bundle over a
  manifold  $X$ of odd dimension $m$. By the {\it twisted
semi-characteristic}
of $X$ (with coefficientes in $F$) we mean the integer
$$s\chi_F(X)= \sum_{i=0}^{(m-1)/2} \dim\, H_{2i}(X;F).$$
In the case of the trivial real line bundle we recover the
semi-characteristic $s\chi (X)\in \Z$ which appeared in Section 4.4.
The next theorem computes $s\chi_F(X) \,(\mod\,2)$ as a function of 
an orthogonal real vector bundle $F$ in the case $m\equiv 1\, (\mod\,4)$.
We refer to  \cite{LMP}, \cite{K},
for other properties  of  the twisted semi-characteristics.

\proclaim{11.2. Theorem}
Let
$X$ be  a  closed connected orientable smooth manifold  of
 dimension $m \equiv 1\, (\mod\,\,4)$.   Let $F $ be a flat    $\R$-vector
bundle  over
 $X$
with orthogonal structure group. Then   $$s\chi_F(X) \,\equiv
\,\<w_1(F)\cup
w_{m-1}(X),[X]\>\, + \,s\chi (X) \cdot \dim \,F\quad  (\mod\,2).$$
\endproclaim

\demo{Proof}  When we add to $F$ the trivial line bundle, both sides
of the formula increase by $s\chi (X)$. Therefore
it is enough to prove the theorem for   even-dimensional $F$.
Set $x=\<w_1(F)\cup
w_{m-1}(X),[X]\> \in \Z/2\Z$.  We should prove that
$ s\chi_F(X)  \equiv  x  \, (\mod\,2)$.

Fix an Euler structure $\xi\in \Eul(X)$  and consider the  torsion
$\tau(X,\xi;F)$, which is an element of $\det\, H_\ast(X; F)$ (see Section 6.1).
By Theorem 7.2 and remarks in Section 4.4,
$$D(\tau(X,\xi; F )) = \tau (X,\xi^\ast; F^\ast )=(\det_{F^\ast}
c(\xi))^{-1}   \tau (X,\xi;
F^\ast )\tag11-1$$
$$= \det_{F} (c(\xi)) \,  \tau (X,\xi;
F^\ast )=(-1)^{x} \tau (X,\xi;
F^\ast ) $$
where $D:\det \, H_\ast(X; F)
\to  \det\, H_\ast(X;F^\ast)$ is the duality
operator
(4-2).

 The flat scalar product on   $F$ gives an
isomorphism of flat vector bundles $\phi: F^\ast \to F$
which induces an isomorphism ${\phi_\ast}:\det \, H_\ast(X; F^\ast)
\to  \det\, H_\ast(X;F)$. By formula
(11-1),
$${\phi_\ast}(D(\tau(X,\xi; F )))=
(-1)^{x} {\phi_\ast}(\tau (X,\xi;
F^\ast ))= (-1)^{x}   \tau (X,\xi;
F  ).$$ Therefore
the number $(-1)^{x}$ equals the degree of the
linear endomorphism ${\phi_\ast} \circ  D$ of the real line  $\det\, H_\ast(X;F)$.
   We  shall  show  below that
the degree of ${\phi_\ast} \circ D$ equals
 $(-1)^{s\chi_F(X)}$.
This   would imply  $s\chi_F(X) = x\, (\mod \, 2) $ and complete the proof of
the
theorem.

We have  $H_\ast(X;F)=\bigoplus_{p }  V_p$, where $V_p= H_p(X;F)\oplus
H_{m-p}(X;F)$ and $p=0, 1, \dots, (m-1)/2$.
We consider each $V_p$ as a graded vector space   with  zero entries
 in degrees $\neq p, m-p$.
 Using the isomorphism
  $\phi: F^\ast \to F$, we can identify $V_p$ with its dual. In this way
 we obtain   duality
operators
$$D_p: \det \, V_p \to \det \, V_p,\quad p=0, 1, \dots, (m-1)/2.$$
Using the fusion isomorphism constructed in Section 2.3 and Lemma
2.7.3 we obtain a natural isomorphism $\mu: \det\, H_\ast(X;F) \to
\otimes_p \, \det V_p$.
Lemma 2.7.2 implies that the conjugation by $\mu$ transforms   ${\phi_\ast} \circ   D$
into the tensor product $ \otimes_p D_p  $. Therefore    $$\deg
({\phi_\ast}   D)=
\prod_{p=0}^{(m-1)/2} \deg D_p.$$  We shall show that   $\deg D_p= (-1)^{\dim
H_p(X;F)}$. This will imply that $\deg ({\phi_\ast}   D)=  (-1)^{s\chi_F(X)}$.

It is clear (from the definitions introduced in Section 2.5) 
that $\deg D_p= (-1)^{s(V_p)}$.  Obviously, we have
$$\alpha_q(V_p) =\cases \dim \, H_p (X;F),\quad\text{if}\quad p\le q
<m-p,\\
0,\quad \text{otherwise}.
\endcases
$$
Thus,  one easily verifies that $\sum_{q=1}^m\alpha_{q-1}(V_p)\alpha_q(V_p)=0 \in \Z/2\Z$
and  
$$
\aligned
&s(V_p)=\sum_{q=0}^{(m-1)/2}\alpha_{2q}(V_p)\\
&={((m+1)}/{2})\cdot \dim \,H_p(X;F)\, \, \, (\mod\, 2)
\endaligned\tag11-2
$$
Therefore we obtain $s(V_p)\equiv \dim \, H_p(X;F)\, (\mod\, 2)$
assuming that $m\equiv 1\, (\mod\,4)$. 
\qed

\enddemo
\subheading{11.3.  Example} Let $X=S^1$. Then  Theorem 11.2
reduces to the  following simple statement:
{\it For $A\in O(n)$,
the codimension of the linear space  of fixed points
of $A$
is even if $A\in SO(n)$, and is odd otherwise.}

\subheading{11.4. Remark} One may use formula (11-2) in the case $m\equiv
3\, (\mod\,4)$
to conclude that $s(V_p)\equiv 0\, \mod \, 2$. Together with the arguments
used in the proof of Theorem 11.2, this gives an
independent proof (which works also for PL manifolds)
of the theorem of Massey \cite{Ma} about vanishing of the Stiefel-Whitney class
$w_{m-1}(X)$ for any closed orientable smooth manifold $X$
of dimension $m\equiv 3\, (\mod\,4)$.


\Refs

\widestnumber\key {BGS}

\ref\key BGS\by J.-M. Bismut, H. Gillet, C. Soul\'e
\paper Analytic torsion and holomorphic determinant bundles, I
\jour Comm. Math. Phys.\vol 115\yr 1988\pages 49-78
\endref

\ref\key BZ\by J.-M. Bismut, W. Zhang\paper
An extension of a theorem by Cheeger and M\"uller
\jour Asterisque\vol 205\yr 1992
\endref

\widestnumber\key {BFK}
\ref\key BFK\by D. Burghelea, L. Friedlander, T. Kappeler\paper
Asymptotic expansion of the Witten deforma\-tion of the analytic torsion
\jour Preprint\yr 1994
\endref

\ref\key C\by J. Cheeger\paper Analytic torsion and the heat equation
\jour Ann. Math. \yr 1979\vol 109\pages 259-322
\endref

\ref\key Fr\by D.S. Freed\paper Reidemeister torsion, spectral sequences,
and Brieskorn spheres\jour J. reine angew. Math. \yr 1992\vol 429
\pages 75-89
\endref

\ref\key Fa\by M. Farber\paper Combinatorial invariants computing the
Ray-Singer
analytic torsion\jour Differential geometry and its applications\vol 6\yr
1996\pages 351-366
\endref

\ref\key HT\by S. Halperin, D. Toledo\paper Stiefel-Whitney homology
classes\jour
Ann. Math. \vol 96\yr 1972\pages 511 - 525\endref

\ref\key K\by G. Kempf\paper Deformations of semi-Euler
characteristics\jour 
Amer. J. Math.
\yr 1992\vol 114\pages 973-978\endref

 \widestnumber\key {LMP}

\ref\key LMP\by G. Luszig, J. Milnor and F. Peterson\paper
Semi-characteristic and
 cobordism
\jour Topology\vol 8\yr 1969\pages 357- 359\endref

\ref\key Ma\by W. Massey\paper On the Shtiefel - Whitney classes of a manifold
\jour Amer. Jour. of Mathematics \vol 82\yr 1960 \pages 92 - 102\endref

\ref\key M1\by J. Milnor\paper A duality theorem for Reidemeister
torsion\jour Ann. of Math.\yr 1962\vol 76\pages 137-147
\endref


\ref\key M2\by J. Milnor\paper Whitehead torsion
\jour Bull. Amer. Math. Soc.
\yr 1966 \vol 72 \pages 358 - 426
\endref

\ref\key Mu1\by W. M\"uller\paper Analytic torsion and R-torsion
for Riemannian manifolds\jour Advances of Math.\yr 1978\vol 28
\pages 233-305
\endref

\ref\key Mu2\by W. M\"uller\paper Analytic torsion and R-torsion 
for unimodular representations
\jour J. Amer. Math. Soc.\vol 6\yr 1993\pages 721-743
\endref

\ref\key RS\by D.B. Ray, I.M. Singer\paper R-torsion and the Laplacian
on Riemannian manifolds\jour Advances in Math.\yr 1971\vol 7\pages 145-210
\endref

\ref\key RS1 \by C.P. Rourke, B.J. Sanderson\book Introduction to
piecewise-linear topology
\publ Springer-Verlag\publaddr Ber\-lin, \-Heidelberg, New York
\yr 1972
\endref


\ref \key T1\by V.G. Turaev\paper Reidemeister torsion in knot theory\jour
Uspekhi Mat. Nauk
41:1(1986), 97-147; English translation: Russian Math. Surveys 41:1(1986),
119-182\endref

\ref \key T2\by V.G. Turaev\paper Euler structures, nonsingular vector
fields, and torsion of
Reidemeister type\jour Izvestia Acad. Sci. USSR 53:3(1989), 130-146;
English translation:
Math. USSR Izvestia 34:3(1990), 627-662\endref

\ref \key T3\by V.G. Turaev\paper Torsion invariants of $Spin^c$-structures
on 3-manifolds
\jour Math. Research Letters \yr 1997\vol 4:5 \pages 679  - 695\endref

\endRefs
\enddocument